\documentclass[11pt]{article} 
\usepackage{latexsym} 
\usepackage{amsmath}
\usepackage{mathptmx}
\usepackage{amssymb} 
\usepackage{amssymb} 
\usepackage{amsfonts}
\usepackage{amscd} 
\usepackage[all]{xy}
\usepackage{graphicx}
\usepackage{layout}
\begin{document} 
%%%%%%%%%%%%%%%%%%%%%%%%%%%%%%%%% 
%%%%%%%%%%%%%%%%%%%%%%%%%%%%%%%%%
\newtheorem{Th}{Theorem}[section]
\newtheorem{Cor}{Corollary}[section]
\newtheorem{Prop}{Proposition}[section]
\newtheorem{Lem}{Lemma}[section]
\newtheorem{Def}{Definition}[section]
\newtheorem{Rem}{Remark}[section]
\newtheorem{Ex}{Example}[section]
\newtheorem{stw}{Proposition}[section]

%Definitions of bet ent

\newcommand{\bet}{\begin{Th}}
\newcommand{\ent}{\stepcounter{Cor}
   \stepcounter{Prop}\stepcounter{Lem}\stepcounter{Def}
   \stepcounter{Rem}\stepcounter{Ex}\end{Th}}

%Definitions of bec enc bep enp bel enl
%bef enf ber enr bee ene

\newcommand{\bec}{\begin{Cor}}
\newcommand{\enc}{\stepcounter{Th}
   \stepcounter{Prop}\stepcounter{Lem}\stepcounter{Def}
   \stepcounter{Rem}\stepcounter{Ex}\end{Cor}}
\newcommand{\bep}{\begin{Prop}}
\newcommand{\enp}{\stepcounter{Th}
   \stepcounter{Cor}\stepcounter{Lem}\stepcounter{Def}
   \stepcounter{Rem}\stepcounter{Ex}\end{Prop}}
\newcommand{\bel}{\begin{Lem}}
\newcommand{\enl}{\stepcounter{Th}
   \stepcounter{Cor}\stepcounter{Prop}\stepcounter{Def}
   \stepcounter{Rem}\stepcounter{Ex}\end{Lem}}
\newcommand{\bef}{\begin{Def}}
\newcommand{\enf}{\stepcounter{Th}
   \stepcounter{Cor}\stepcounter{Prop}\stepcounter{Lem}
   \stepcounter{Rem}\stepcounter{Ex}\end{Def}}
\newcommand{\ber}{\begin{Rem}}
\newcommand{\enr}{
   %\stepcounter{Rem} 
   \stepcounter{Th}\stepcounter{Cor}\stepcounter{Prop}
   \stepcounter{Lem}\stepcounter{Def}\stepcounter{Ex}\end{Rem}}
\newcommand{\bee}{\begin{Ex}}
\newcommand{\ene}{
 %\stepcounter{Ex}
   \stepcounter{Th}\stepcounter{Cor}\stepcounter{Prop}
   \stepcounter{Lem}\stepcounter{Def}\stepcounter{Rem}\end{Ex}}
\newcommand{\Proof}{\noindent{\it Proof\,}:\ }
\newcommand{\beP}{\Proof}
\newcommand{\enP}{\hfill $\Box$ \par\vspace{5truemm}}
%\newcommand{\enP}{\QED}

%%%%%%%%%%%%%%%%%%%%%%%%%%%%%%%%%%%%%%%%%
%Beginning of Local Definition
%Local definitions
\newcommand{\EE}{\mathbb{E}}
\newcommand{\QQ}{\mathbb{Q}}
\newcommand{\R}{\mathbb{R}}
\newcommand{\C}{\mathbb{C}}
\newcommand{\ZZ}{\mathbb{Z}}
\newcommand{\KK}{\mathbb{K}}
\newcommand{\NN}{\mathbb{N}}
\newcommand{\PP}{\mathbb{P}}
\newcommand{\HH}{\mathbb{H}}
\newcommand{\uuu}{\boldsymbol{u}}
\newcommand{\xxx}{\boldsymbol{x}}
\newcommand{\aaa}{\boldsymbol{a}}
\newcommand{\bbb}{\boldsymbol{b}}
\newcommand{\AAA}{\mathbf{A}}
\newcommand{\BBB}{\mathbf{B}}
\newcommand{\LLL}{\mathbf{L}}
\newcommand{\ccc}{\boldsymbol{c}}
\newcommand{\iii}{\boldsymbol{i}}
\newcommand{\jjj}{\boldsymbol{j}}
\newcommand{\kkk}{\boldsymbol{k}}
\newcommand{\rrr}{\boldsymbol{r}}
\newcommand{\FFF}{\boldsymbol{F}}
\newcommand{\yyy}{\boldsymbol{y}}
\newcommand{\ppp}{\boldsymbol{p}}
\newcommand{\qqq}{\boldsymbol{q}}
\newcommand{\nnn}{\boldsymbol{n}}
\newcommand{\vvv}{\boldsymbol{v}}
\newcommand{\eee}{\boldsymbol{e}}
\newcommand{\fff}{\boldsymbol{f}}
\newcommand{\www}{\boldsymbol{w}}
\newcommand{\0}{\boldsymbol{0}}
\newcommand{\lon}{\longrightarrow}
\newcommand{\ga}{\gamma}
\newcommand{\pa}{\partial}
\newcommand{\QED}{\hfill $\Box$}
\newcommand{\id}{{\mbox {\rm id}}}
\newcommand{\Ker}{{\mbox {\rm Ker}}}
\newcommand{\grad}{{\mbox {\rm grad}}}
\newcommand{\ind}{{\mbox {\rm ind}}}
\newcommand{\rot}{{\mbox {\rm rot}}}
\newcommand{\diver}{{\mbox {\rm div}}}
\newcommand{\Gr}{{\mbox {\rm Gr}}}
\newcommand{\GL}{{\mbox {\rm GL}}}
\newcommand{\LG}{{\mbox {\rm LG}}}
\newcommand{\Diff}{{\mbox {\rm Diff}}}
\newcommand{\Symp}{{\mbox {\rm Symp}}}
\newcommand{\Ct}{{\mbox {\rm Ct}}}
\newcommand{\Uns}{{\mbox {\rm Uns}}}
\newcommand{\rank}{{\mbox {\rm rank}}}
\newcommand{\sign}{{\mbox {\rm sign}}}
\newcommand{\Spin}{{\mbox {\rm Spin}}}
\newcommand{\Sp}{{\mbox {\rm Sp}}}
\newcommand{\Int}{{\mbox {\rm Int}}}
\newcommand{\Hom}{{\mbox {\rm Hom}}}
\newcommand{\Tan}{{\mbox {\rm Tan}}}
\newcommand{\codim}{{\mbox {\rm codim}}}
\newcommand{\ord}{{\mbox {\rm ord}}}
\newcommand{\Iso}{{\mbox {\rm Iso}}}
\newcommand{\corank}{{\mbox {\rm corank}}}
\def\mod{{\mbox {\rm mod}}}
\newcommand{\pt}{{\mbox {\rm pt}}}
\newcommand{\qed}{\hfill $\Box$ \par}
\newcommand{\spe}{\vspace{0.4truecm}}
\renewcommand{\0}{\mathbf 0}
\newcommand{\ad}{{\mbox{\rm ad}}}
\newcommand{\xdownarrow}[1]{%
  {\left\downarrow\vbox to #1{}\right.\kern-\nulldelimiterspace}
}

\newcommand{\dint}[2]{{\displaystyle\int}_{{\hspace{-1.9truemm}}{#1}}^{#2}}
%%%%%%%%%%%%%%%%%%%%%%%%%%%%%%%%%%%%%%%%%%%%%%%%%
%End of local definitions

%%%TITLE
\title{Recognition Problem of Frontal Singularities
}

\author{Goo \textsc{Ishikawa}
\thanks{Faculty of Science, Hokkaido University, Sapporo 060-0810, Japan.
ishikawa@math.sci.hokudai.ac.jp
}}

%\address{Department of Mathematics, Hokkaido University, Sapporo 060-0810, Japan. }
%\email{ishikawa@math.sci.hokudai.ac.jp}

%\subjclass[2010]{Primary 57R45; Secondary 58K50, 53A07, 53D12, 53C50}% Subject code(s)
%
%\keywords{Jacobi ideal, Jacobi module, ramification module, kernel field, opening, Lorentzian manifold}% Key word(s)

\renewcommand{\thefootnote}{\fnsymbol{footnote}}
\footnotetext{
\\
\qquad
The author was supported by JSPS KAKENHI No.15H03615 and No.15K13431.}

\date{}

\maketitle

\section{Introduction}

This is a survey article on recognition problem of frontal singularities. 

First we explain the recognition problem of singularities and 
its significance. 

Let $f : (\R^n, a) \to (\R^m, b)$ and $f' : (\R^n, a') \to (\R^m, b')$ 
be smooth (= $C^\infty$) map-germs. 
Then $f$ and $f'$ are called {\it ${\mathcal A}$-equivalent}
or {\it diffeomorphic} if there exist diffeomorphism-germs 
$\sigma : (\R^n, a) \to (\R^n, a')$ and $\tau : (\R^m, b) \to (\R^m, b')$ such that 
the diagram 
$$
\begin{array}{ccc}
(\R^n, a) & \xrightarrow{\ f\ } & (\R^m, b) 
%\vspace{0.4truecm}
\\
\ \ \downarrow \sigma & & \ \ \downarrow \tau
%\vspace{0.2truecm}
\\
(\R^n, a') & \xrightarrow{\  f'\ } & (\R^m, b') 
\end{array}
$$
commutes. 

By a {\it singularity} of smooth mappings, 
we mean an {\it ${\mathcal A}$-equivalence class} of map-germs. 

\

Suppose that we investigate \lq\lq singularities\rq\rq of mappings belonging to some given class. 
Then the recognition problem of singularities may be understood as the following dual manners: 

\

\noindent
{\it Problem: }
Given two map-germs $f$ and $f'$, {\it belonging to the given class, }
determine, {\it as easily as possible} 
whether $f$ and $f'$ are equivalent or not. 

\noindent
{\it Problem: }
Given a singularity, find criteria, 
{\it as easy as possible}, 
to determine whether a map-germ $f$ {\it belonging to some class} 
has (= falls into) the given singularity or not. 

\

Importance of the recognition problem of singularities can be explained as follows. 

Once we establish a classification list of singularities in a {\it situation A}, 
we will face (at least) {\it two kinds of needs}: 

1. Given a map-germ in the same {\it situation A}, we want to know which singularity is it 
in the list. 

2. For another {\it situation B}, we want to know how similar is the classification 
list of singularities as {\it A} or not. 

In both cases, we need to recognize the singularities, {\it as easily as possible}, by {\it as many as possible} 
criteria. 
It is indispensable the recognition of singularities for applications of singularity theory, and 
to solve classification problems in various situations. 

\

The recognition problem of singularities of smooth map-germs has been 
treated by the many mathematicians, 
motivated by differential geometry and other wide area, 
and its solutions are supposed to have many applications. 

In fact most of known results of recognition of singularities are found under the motivation of 
geometric studies of singularities appearing in Euclid geometry and various Klein geometries 
(\cite{KRSUY}\cite{FSUY}\cite{IST}). 

\bee
{\rm 
(Singularities in non-Euclidean geometry) \ 
The following is a diagram representing the history of non-Euclidean geometry 
found in the reference \cite{Sharp}: 
%\begin{center}
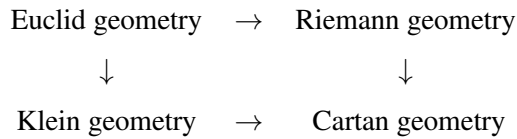
$$
\begin{tabular}{ccc}
{\mbox{\rm Euclid geometry}} & {\mbox{\rm $\rightarrow$}} & {\mbox{\rm Riemann geometry}} 
\vspace{0.2truecm}
\\
{\mbox{\rm $\downarrow$}} & & {\mbox{\rm $\downarrow$}}
\vspace{0.2truecm}
\\
{\mbox{\rm Klein geometry}} & {\mbox{\rm $\rightarrow$}} & {\mbox{\rm Cartan geometry}} 
\end{tabular}
$$
Then it would be natural to ask 

\

\noindent
{\it Problem:} 
How are the classification results of singularities in Euclid geometry (resp. in Klein geometry) 
valid in Riemann geometry (resp. in Cartan geometry)? 

\noindent
In other words, 
\\
{\it Problem:} 
Do the classifications of singularities in flat ambient spaces work also for \lq\lq curved" ambient spaces? 

\

%To solve the classification problem it is indispensable the recognition of singularities. 
%In fact we applied the several results of recognition (\cite{KRSUY}\cite{FSUY}) to, 
%for instance, the classification problem of generic singularities appearing in tangent surfaces which are 
%ruled by geodesics in general Riemannian spaces (\cite{IY17}\cite{IY18}). 
In fact, we applied the several results of recognition(\cite{KRSUY}\cite{FSUY}), 
for instance, to the generic classification of singularities of improper affine spheres 
and of surfaces of constant Gaussian curvature(\cite{IM}), and moreover, to 
the classification of generic singularities appearing in tangent surfaces which are 
ruled by geodesics in general Riemannian spaces (\cite{IY17}\cite{IY18}). See also \S\ref{application}. 
}
\ene

In this paper we will pay our attention to the class of mappings, {\it frontal mappings}, which is introduced and 
studied in \S\ref{frontals}. Then we survey several 
recognition theorems on them in \S\ref{Recognition}. 
Note that the recognitions of fronts or frontals $(\R^n, a) \to \R^m$ 
are studied by many authors (\cite{KRSUY}\cite{FSUY}\cite{Saji}\cite{SUY}\cite{Kabata}). 

To show the theorems given in \S\ref{Recognition}, 
we introduce the notion of {\it openings}, relating it with that of frontals, in \S\ref{Frontals and openings}. 
See also \cite{Ishikawa13}\cite{Ishikawa14}. 
In fact, in \S\ref{Frontals and openings}, 
we observe that {\it any frontal singularity is an opening of a map-germ from 
$\R^n$ to $\R^n$} (Lemma \ref{opening-frontal}). 

Then we naturally propose: 

\

\noindent
{\it Problem:} 
Study the {\it recognition problem} of frontals from the 
{\it recognition results} on map-germs $(\R^n, a) \to \R^n$, ({\it $n = m$}), 
combined with the viewpoint of {\it openings}. 

\

In this paper, in connection with the above problems, 
we specify {\it geometrically} several frontal singularities which we are going to 
treat (Example \ref{tangent-surfaces}). 
Then we solve the recognition problem of such singularities, in \S\ref{Recognition}, 
giving explicit normal forms. 
In fact we combine the recognition results on $(\R^2, 0) \to (\R^2, 0)$ by K. Saji ($\sim$2010) 
and several arguments on openings, which was 
implicitly performed for the classification of singularities of tangent surfaces 
(tangent developables) by the author ($\sim$1995) over twenty years, 
the idea of which traces back to the author's master thesis \cite{Ishikawa83}. 
We prove recognition theorems in \S\ref{Proofs-of-recognition-theorems}. 

In the last section \S\ref{application}, 
as an application of our solutions of recognition problem of frontal singularities, 
we announce the classification 
of singularities appearing in tangent surfaces of generic null curves which are ruled by null geodesics 
in general Lorentz $3$-manifolds (\cite{IMT11}\cite{IMT18}), mentioning related recognition results 
and open problems. 

\

In this paper, all manifolds and mappings are assumed to be of class $C^\infty$ unless otherwise stated. 

The author deeply thanks to an anonymous referee for his/her helpful comments to improve the paper.

\section{Frontal singularities}
\label{frontals}

Let $f : (\R^n, a) \to (\R^m, b)$ be a map-germ. Suppose $n \leq m$. 

Then $f$ is called a {\it frontal map-germ} or a {\it frontal} in short,  if 
there exists a smooth ($C^\infty$) family of $n$-planes ${\widetilde f}(t) \subseteq T_{f(t)}\R^m$ 
along $f$, $t \in (\R^n, a)$, i.e. 
there exists a smooth lift ${\widetilde f} : (\R^n, a) \to {\mbox{\rm Gr}}(n, T\R^m)$ 
satisfying the \lq\lq {\it integrality condition}\rq\rq\,  
$$
T_tf(T_t\R^n) \subset \ {\widetilde f}(t) \ (\subset T_{f(t)}\R^m), 
$$
for any $t \in \R^n$ nearby $a$,  such that $\pi\circ \widetilde{f} = f$: 
$$
\xymatrix{
 & {\mbox{\rm Gr}}(n, T\R^m) \ar[d]^{\pi}
 \\
(\R^n, a) \ar[r]_f \ar[ru]^{\widetilde{f}}  & (\R^m, b). 
%    G \ar[r]^{f} \ar[d]_{\pi} & G' \\
%    G/ker f \ar@{.>}[ru]_{\psi}
}
$$
Here ${\mbox{\rm Gr}}(n, T\R^m)$ is the {\it Grassmann bundle} consisting 
of $n$-planes $V \subset T_x\R^m, (x \in \R^m)$ with the canonical projection 
$\pi(x, V) = x$, and $T_tf : T_t\R^n \to T_{f(t)}\R^m$ is the differential of $f$ at $t \in (\R^n, a)$. 

Then ${\widetilde{f}}$ is called a {\it Legendre lift} or an {\it integral lift} of the frontal $f$. 
Actually ${\widetilde{f}}$ is an integral mapping to the canonical or contact distribution on 
${\mbox{\rm Gr}}(n, T\R^m)$ (cf. \cite{Ishikawa12}).

\bee
\label{example-of-frontals}
{\rm 
(1) Any {\it immersion} is a frontal.   In fact then the Legendre lift is given by $\widetilde{f}(t) := T_tf(T_t\R^n)$. 

(2) Any map-germ $(\R^n, a) \to (\R^n, b), ({\it n = m})$ is a frontal. 
In fact the Legendre lift is given by $\widetilde{f}(t) := T_{f(t)}\R^n$.  

(3) Any {\it constant map-germ} is a frontal. In fact we can take any lift $\widetilde{f}$ of $f$.

(4) Any {\it wave-front} $(\R^n, a) \to (\R^{n+1}, b)$, that is a Legendre projection 
of a Legendre submanifold in $\Gr(n, T\R^{n+1}) = PT^*\R^{n+1}$, is a frontal. Take  
the inclusion of the Legendre submanifold as the Legendre lift. 
}
\ene

\bee

{\rm 
\label{tangent-surfaces}
(Singularities of tangent surfaces) \ 
Let $\gamma : (\R, 0) \to \R^m$ be a curve-germ in Euclidean space. 
Then the tangent surface 
${\mbox{\rm Tan}}(\gamma) : (\R^2, 0) \to \R^m$ is defined as the ruled surface 
generated by tangent lines along the curve. 
Suppose $\gamma$ is of type $\LLL = (\ell_1, \ell_2, \ell_3, \dots, ), (1 \leq \ell_1 < \ell_2 < \ell_3 < \cdots)$, i.e. 
$$
\gamma(t) = (t^{\ell_1} + \cdots, \ t^{\ell_2} + \cdots, \ t^{\ell_3} + \cdots, \ \dots ) 
$$
for an affine coordinates of $\R^m$ centered at $\gamma(0)$. Then it is known that 
the singularity of ${\mbox{\rm Tan}}(\gamma)$ 
is {\it uniquely determined} by the type $\LLL$ and called cuspidal edge ({\it CE}) if $\LLL = (1, 2, 3, \dots)$, 
folded umbrella ({\it FU}) or cuspidal cross cap ({\it CCC}) if $(1, 2, 4)$, 
swallowtail ({\it SW}) if $(2, 3, 4)$, Mond ({\it MD}) or cuspidal beaks ({\it CB}) if $(1, 3, 4)$, 
Shcherbak ({\it SB}) if $(1, 3, 5)$, 
cuspidal swallowtail ({\it CS}) if $(3, 4, 5)$, 
open folded umbrella ({\it OFU}) if $(1, 2, 4, 5, \dots)$, open swallowtail ({\it OSW}) if 
$(2, 3, 4, 5, \dots)$, open Mond ({\it OMD}) or open cuspidal beaks ({\it OCB}) if $(1, 3, 4, 5, \dots)$ (see \cite{Ishikawa12}). 
}
\ene

\

In general, a frontal $f : (\R^n, a) \to (\R^m, b)$ is called a {\it front} if 
$f$ has an immersive Legendre lift $\widetilde{f}$. 

\

Let ${\mathcal E}_a := \{ h : (\R^n, a) \to \R\}$ 
denote the $\R$-algebra of smooth function-germs 
on $(\R^n, a)$. 

Denote by $\Gamma$ the set of subsets $I \subseteq \{ 1, 2, \dots, m\}$ 
with $\#(I) = n$. 
For a map-germ $f : (\R^n, a) \to (\R^m, b), n \leq m$ and $I \in \Gamma$, 
we set $D_I = \det(\pa f_i/\pa t_j)_{i \in I, 1 \leq j \leq n}$. 
Then {\it Jacobi ideal} $J_f$ of $f$ is defined as the ideal 
generated in ${\mathcal E}_a$ by all $n$-minor determinants $D_I$ ($I \in \Gamma$) of Jacobi matrix $J(f)$ of $f$. 
Then we have: 

\bel
\label{Criterion-of-frontality}
{\rm (Criterion of frontality)} 
Let $f : (\R^n, a) \to (\R^m, b)$ be a map-germ. 
If $f$ is a frontal, then the Jacobi ideal $J_f$ of $f$ is principal, 
i.e. it is generated by one element. In fact $J_f$ is generated by $D_I$ for some $I \in \Gamma$. 
Conversely, if $J_f$ is principal and the singular locus 
$$
S(f) = \{ t \in (\R^n, a) \mid \rank(T_tf : T_t\R^n \to T_{f(t)}\R^m) < n \}
$$ 
of $f$ is nowhere dense in $(\R^n, a)$, then $f$ is a frontal. 
\enl

%First we give the proof of Lemma \ref{Criterion-of-frontality}. 
%
%\noindent
%{\it Proof of Lemma \ref{Criterion-of-frontality}:} 
\Proof
Let $f$ be a frontal and $\widetilde{f}$ be a Legendre lift of $f$. 
Take $I_0 \in \Gamma$ such that $\widetilde{f}(a)$ 
projects isomorphically by the projection 
$\R^m \to \R^n$ to the components belonging to $I_0$. Let 
$(p_I)_{I \in \Gamma}$ 
be the Pl{\" u}cker coordinates of $\widetilde{f}$. Then $p_{I_0}(a) \not= 0$. 
This implies that for any $I \in \Gamma$, there exists $h_I \in {\mathcal E}_a$ such that 
$D_I = h_I D_{I_0}$. Set $\lambda = D_{I_0}$. Then the Jacobi ideal $J_f$ 
is generated by $\lambda$. 

Conversely suppose $J_f$ is generated by one element $\lambda \in {\mathcal E}_a$. 
Since $J_f$ is generated by $\lambda$, we have that there exists $k_I \in {\mathcal E}_a$ for any $I \in \Gamma$ 
such that $D_I = k_I \lambda$. 
Since $\lambda \in J_f$, 
there exists $\ell_I \in {\mathcal E}_a$ for any $I \in \Gamma$ 
such that $\lambda = \sum_{I \in \Gamma} \ell_I D_I$. 
Therefore $(1 - \sum_{I \in \Gamma} \ell_Ik_I)\lambda = 0$. 
Suppose $(\ell_Ik_I)(a) = 0$ for any $I \in \Gamma$. Then $1 - \sum_{I \in \Gamma} \ell_Ik_I$ is 
a unit and therefore $\lambda = 0$. Thus we have $J_f = 0$. This contradicts to the assumption 
that $S(f)$ is nowhere dense. 
Hence there exists $I_0 \in \Gamma$ such that $(\ell_{I_0}k_{I_0})(a) \not= 0$. 
Then $k_{I_0}(a) \not= 0$. Therefore $J_f$ is generated by $D_{I_0}$. Hence 
$D_I = h_I D_{I_0}$ for any $I \in \Gamma$ with $h_{I_0}(a) = 1$. 
Then the Legendre lift $\widetilde{f}$ on $\R^n \setminus S(f)$ 
extends to $(\R^n, a)$, which is given by the Pl{\" u}cker coordinates $(h_I)_{I \in \Gamma}$. 
\QED

\bee
{\rm 
Define $f : (\R^2, 0) \to (\R^3, 0)$ by $f(t_1, t_2) := (\varphi(t_1), \ \varphi(t_1)t_2, \ \varphi(-t_1))$, 
where the $C^\infty$ function 
$\varphi : (\R, 0) \to (\R, 0)$ is given by $\varphi(t) = \exp(-1/t^2) (t \geq 0), 0 (t \leq 0)$. 
Then the Jacobi ideal $J_f$ is generated by $\varphi'(t_1)\varphi(t_1)$ and therefore $J_f$ is principal and $J_f \not= 0$. 
However $f$ is not a frontal. In fact, for $t_1 > 0$, $(T_{(t_1,t_2)}f)(T_{(t_1, t_2)}\R^2)$ is given by the plane $dx_3 = 0$ and 
for $t_1 < 0$, $(T_{(t_1,t_2)}f)(T_{(t_1, t_2)}\R^2)$ contains the $x_3$-axis. Therefore $f$ can not be a frontal. 
}
\ene

\bec
Let $f : (\R^n, a) \to (\R^m, b)$ be a map-germ. Suppose $f$ is analytic and $J_f \not= 0$. 
Then $f$ is a frontal if and only if $J_f$ is a principal ideal. 
\enc

\Proof
By Lemma \ref{Criterion-of-frontality}, if $f$ is frontal, then $J_f$ is principal. 
If $J_f$ is principal and $J_f \not= 0$, then $D_I \not= 0$ for some $I \in \Gamma$. 
Since $f$ is analytic, $S(f)$ is nowhere dense. Thus by Lemma \ref{Criterion-of-frontality}, 
$f$ is a frontal. 
\QED

\bee
{\rm 
Define $f : (\R^3, 0) \to (\R^4, 0)$ by $f(t_1, t_2, t_3) := (t_1^3, \ t_1^2t_2, \ t_1t_2^2, \ t_2^3)$. 
The germ $f$ parametrizes the cone over a non-degenerate cubic in $P(\R^4) = \R P^3$. 
Then $f$ is analytic and $J_f = 0$ is principal. However $f$ is not a frontal. 
}
\ene

\bef
\label{Jacobian}
{\rm 
Let $f : (\R^n, a) \to (\R^m, b)$ be a frontal. 
Then a generator $\lambda \in {\mathcal E}_a$ 
of $J_f$ is called a {\it Jacobian} (or a {\it singularity identifier}) of $f$, 
which is uniquely determined from $f$ up to multiplication of a unit in ${\mathcal E}_a$. 
}
\enf

The singular locus $S(f)$ of a frontal $f$ is given by the zero-locus of the Jacobian $\lambda$ of $f$. 

\bef 
{\rm (Proper frontals) }
{\rm 
A frontal $f : (\R^n, a) \to (\R^m, b)$ is called {\it proper} if 
the singular locus 
$S(f)$ is nowhere dense in $(\R^n, a)$. 
}
\enf

\ber
{\rm 
Our naming \lq\lq proper\rq\rq\, is a little confusing since its usage is different from 
the ordinary meaning of properness (inverse images of any compact is compact). 
Our condition that the singular locus $S_f$ is nowhere dense is easy to handle for the local study of mappings. 
}
\enr

\bel
Let $f : (\R^n, a) \to (\R^m, b)$ be a proper frontal or $n = m$. 
Then $f$ has a unique Legendre lift $\widetilde{f} :  (\R^n, a) \to \Gr(n, T\R^m)$. 
\enl

\Proof
On the regular locus $\R^n \setminus S(f)$, there is 
the unique Legendre lift $\widetilde{f}$ defined by 
$\widetilde{f}(t) := (T_tf)(T_t\R^n)$. 
Let $f$ be a proper frontal. Then $\R^n \setminus S(f)$ 
is dense in $(\R^n, a)$.  Therefore the extension of $\widetilde{f}(t)$ is unique. 
Let $n = m$. Then the unique lift $\widetilde{f}$ is defined by $\widetilde{f}(t) = T_{f(t)}\R^m$ 
(Example \ref{example-of-frontals} (2)). 
\QED

\

Let $f : (\R^n, a) \to (\R^m, b)$ be a frontal
(resp. a proper frontal) 
and $\widetilde{f} : (\R^n, a) \to {\mbox{\rm Gr}}(n, T\R^m)$ a Legendre lift of $f$. 
Recall that $\widetilde{f}(t), (t \in (\R^n, a))$ is an 
$n$-plane field along $f$. In particular $\widetilde{f}(a) \subseteq T_b\R^m$. 

\bef
\label{adapted-coordinates}
{\rm 
A system $(x_1, \dots, x_n, x_{n+1}, \dots, x_m)$ of local coordinates of $\R^m$ centered at 
$b$ is called {\it adapted} to 
$\widetilde{f}$ (or, to $f$) 
 if 
$$
\begin{array}{rcl}
\widetilde{f}(a) & = & \left\langle \left(\dfrac{\pa}{\pa x_1}\right)_b, 
\dots, \left(\dfrac{\pa}{\pa x_n}\right)_b\right\rangle_{\R}
\vspace{1truecm}
\\
& ( = & \{ v \in T_b\R^m \mid dx_{n+1}(v) = 0, \dots, dx_{m}(v) = 0 \} ). 
\end{array}
$$
}
\enf

Clearly we have 

\bel
Any frontal $f : (\R^n, a) \to (\R^m, b)$ 
has an adapted system of local coordinates on $(\R^m, b)$. In fact any 
system of local coordinates on $(\R^m, b)$ is modified into an adapted system of local coordinates 
by a linear change of coordinates. 
\enl

\ber
{\rm 
For an adapted system of coordinates $(x_1, \dots, x_n, x_{n+1}, \dots, x_m)$ of $f$, 
the Jacobian $\lambda$ is given by the ordinary Jacobian $\frac{\pa(f_1, \dots, f_n)}{\pa(t_1, \dots, t_n)}$, 
where $f_i = x_i\circ f$. 
}
\enr

\bee
{\rm 
Let $f : (\R^2, 0) \to (\R^3, 0)$ be given by 
$$
(u, t) \mapsto (x_1, x_2, x_3) = (t + u, \ t^3 + 3t^2u, \ t^4 + 4t^3u), 
$$
which is the tangent surface, Mond surface, of the curve $t \mapsto (t, t^3, t^4)$. 

Then the Jacobi matrix $J(f)$ of $f$ is given by 
$$
J(f) = \left(
\begin{array}{cc}
1 & 1 
\\
3t^2 & 3t^2 + 6tu 
\\
4t^3 & 4t^3 + 12t^2u
\end{array}
\right), 
\vspace{-0.3truecm}
$$
and its minors are calculated as 
$$
\left\{ 
\begin{array}{rcl}
D_{12} & = & 6tu, 
\vspace{0.2truecm}
\\
D_{13} & = & 12t^2u = 2t(6tu), 
\vspace{0.2truecm}
\\
D_{23} & = & 12t^4u = 2t^3(6tu), 
\end{array}
\right.
$$ 
Then the Jacobi ideal $J_f$ is generated by $\lambda = tu$. 
Therefore $f$ is a proper frontal with $S(f) = \{ (u, t) \mid tu = 0\}$. 
The unique Legendre lift $\widetilde{f} : (\R^2, 0) 
\to {\mbox{\rm Gr}}(2, T\R^3)$ of $f$ is given, via the Pl{\" u}cker coordinates of fibre components, 
$$
D_{12}/D_{12} = 1, \ D_{13}/D_{12} = 2t, \ D_{23}/D_{12} = 2t^3. 
$$
The system of coordinates $(x_1, x_2, x_3)$ is adapted for $f$ in the example. 
}
\ene

\section{Recognition of several frontal singularities}
\label{Recognition}

To give our recognition results we need the notion of \lq\lq kernel fields" in addition to that of Jacobians 
of frontals. 

\

Let $f : (\R^n, a) \to (\R^m, b)$ be a map-germ. We denote by ${\mathcal V}_a$ the ${\mathcal E}_a$-module 
of vector fields over $(\R^n, a)$ and set 
$$
{\mathcal N}_f := \{ \eta \in {\mathcal V}_a \mid \eta f_i \in J_f, \ (1 \leq i \leq m)\}, 
$$
which is a ${\mathcal E}_a$-submodule of ${\mathcal V}_a$. 

Note that, if $\eta \in {\mathcal {N}}_f$,  
then $\eta(t) \in \Ker(T_tf : T_t \R^n \to T_{f(t)}\R^m)$ for any $t \in S(f)$. 
Moreover note that, if $\lambda \in J_f$, then $\lambda\cdot {\mathcal V}_a \subseteq {\mathcal N}_f$. 

A map-germ $f : (\R^n, a) \to (\R^m, b)$ is called {\it of corank} $k$ if 
$\dim_{\R}\Ker(T_a f : T_a\R^n \to T_b\R^m) = k$. 

Then we have 

\bel
\label{kernel-field-free}
Let $f : (\R^n, a) \to (\R^m, b)$ be a map-germ of corank $1$. 
Then ${\mathcal N}_f/J_f\cdot{\mathcal V}_a$ is a free ${\mathcal E}_a$-module of rank $1$, 
i.e. 
${\mathcal N}_f/J_f\cdot{\mathcal V}_a$ is isomorphic to ${\mathcal E}_a$ as 
${\mathcal E}_a$-modules by $[\eta] \to 1$, for some $\eta \in {\mathcal N}_f$. 
\enl

Let $f : (\R^n, a) \to (\R^m, b)$ be a frontal of corank $1$ 
and $\lambda_f$ the Jacobian of $f$ (Definition \ref{Jacobian}). 
Then by Lemma \ref{kernel-field-free}, 
${\mathcal N}_f/\lambda_f\cdot{\mathcal V}_a$ is a free module of rank $1$. 

\bef
{\rm 
A vector field $\eta$ over $(\R^n, a)$ is called a {\it kernel field} (or a {\it null field}) of $f$ if 
$\eta$ generates the free ${\mathcal E}_a$-module ${\mathcal N}_f/\lambda_f\cdot{\mathcal V}_a$. 
}
\enf

\ber
{\rm 
The notion of null fields is introduced first in \cite{KRSUY}. 
}
\enr

\

\noindent
{\it Proof of Lemma \ref{kernel-field-free}:} 
Since $f$ is of corank $1$, $f$ is ${\mathcal A}$-equivalent to a map-germ 
$(\R^n, 0) \to (\R^m, 0)$ of form 
$$
g = (t_1, \dots, t_{n-1}, \varphi_n(t), \dots, \varphi_m(t)). 
$$
Note that ${\mathcal N}_f/J_f{\mathcal V}_a$ 
is isomorphic to ${\mathcal N}_g/J_f{\mathcal V}_0$. 
Moreover 
the Jacobian ideal of $g$ is generated by  
$$
\pa\varphi_n(t)/\pa t_n, \dots, \pa\varphi_m(t)/\pa t_n. 
$$
Let $\eta = \sum_{i=1}^n \eta_i \pa/\pa t_i \in {\mathcal V}_0$. Then $\eta \in 
{\mathcal N}_g$ if and only if $\eta_1, \dots, \eta_{n-1} \in J_g$. 
Therefore ${\mathcal N}_g/J_g{\mathcal V}_0$ is freely generated by 
$\pa/\pa t_n$. Thus we have that ${\mathcal N}_f/J_f\cdot{\mathcal V}_a$ 
is a free ${\mathcal E}_a$-module of rank $1$, 
\QED

\

Now we start to give our recognition theorems on the frontal singularities introduced in 
Example \ref{tangent-surfaces}. 
To begin with, we recall the following fundamental recognition result due to Saji (\cite{Saji}), 
which is a reformulation of Whitney's original results in \cite{Whitney} for parts (1) and (2). 

\bet
\label{Whitney, Saji}
{\rm (Saji\cite{Saji})}
Let $f : (\R^2, a) \to (\R^2, b)$ be a (frontal) map-germ of corank $1$. 
Then, for the Jacobian $\lambda$ and the kernel field $\eta$ of $f$, we have

{\rm (1)} $f$ is ${\mathcal A}$-equivalent to the fold, i.e. to $(t_1, t_2) \mapsto (t_1, t_2^2)$,  
if and only if $(\eta\lambda)(a) \not= 0$. 

{\rm (2)} $f$ is ${\mathcal A}$-equivalent to Whitney's cusp, i.e. to 
$(t_1, t_2) \mapsto (t_1, t_2^3 + t_1t_2)$, if and only if 
$(d\lambda)(a) \not= 0, (\eta\lambda)(a) = 0, (\eta\eta\lambda)(a) \not= 0$. 

{\rm (3)} $f$ is ${\mathcal A}$-equivalent to bec {\` a} bec (beak-to-beak), 
$(t_1, t_2) \mapsto (t_1, t_2^3 + t_1t_2^2)$, 
if and only if $\lambda$ has an indefinite Morse critical point at $a$ 
and $(\eta\eta\lambda)(a) \not= 0$. 
\ent

%\ber
%{\rm 
%The name bec {\` a} bec, that is due to Run{\' e} Thom, 
%can be translated into Japanese "seppun", which also means beak-to-beak. 
%}
%\enr
%

\ber
{\rm 
Each condition (1), (2), (3) of Theorem \ref{Whitney, Saji} is independent of the choice 
of $\lambda$ and $\eta$ and depends only on ${\mathcal J}$-equivalence class of $f$ 
which is introduced in Definition \ref{J-equivalence}. 
In fact, if ${\mathcal J}_{f'\circ\sigma} = {\mathcal J}_f$, then $f'$ satisfies the condition 
for $\lambda' = \lambda\circ\sigma^{-1}$ and $\eta' = (T\sigma)\eta\circ\sigma^{-1}$. (See \S\ref{Frontals and openings}). 
}
\enr

\ber
{\rm 
For a map-germ $f : (\R^2, a) \to (\R^2, b)$ of corank $1$, 
the condition $(d\lambda)(a) \not= 0$ is equivalent to that the Jacobian is ${\mathcal K}$-equivalent 
to the germ $(t_1, t_2) \mapsto t_1$ at the origin. 
The condition that $\lambda$ has an indefinite Morse critical point at $a$ is equivalent to 
that $\lambda$ is ${\mathcal K}$-equivalent 
to the germ $(t_1, t_2) \mapsto t_1t_2$ at the origin. 
}
\enr

\ber
{\rm 
For plane to plane map-germs, the fold (resp. Whitney cusp, bec {\` a} bec) is 
characterized as a\lq\lq tangent map" of a planar curve of type $(1, 2)$ (resp. $(2, 3)$, $(1, 3)$), 
which is ruled by tangent lines to the curve (\cite{Ishikawa12}\cite{IMT16}). 
}
\enr

Let $f : (\R^2, a) \to (\R^m, b), (m \geq 3)$ be a 
proper frontal of corank $1$. We wish to recognize the singularity, i.e. 
${\mathcal A}$-equivalence class of $f$ 
%by its ${\mathcal J}$-equivalence class, 
%in other words, 
by the Jacobian $\lambda = \lambda_f$ and the 
kernel field $\eta = \eta_f$. 
Moreover we wish to recognize the singularity of $f$ as 
an opening of a plane-to-plane map-germ. 
To realize this, we will 
use an adapted system of coordinates $(x_1, x_2, x_3, \dots, x_m)$ for $f$ 
and set $f_i = x_i\circ f$. 
Note that we mention several conditions to recognize singularities 
in terms of adapted coordinates, 
however the conditions are, of course, independent of the choice of an adapted coordinates, 
and therefore any system of adapted coordinates can be taken 
to simplify the checking of a suitable condition. 

%To proceed to the next theorem we introduce the following notion: 

In general, we use the following notation: 

\bef
\label{vanishing-order}
{\rm 
For a germ of vector field $\eta \in {\mathcal V}_a$ over $(\R^n, a)$ and a function-germ 
$h \in {\mathcal E}_a$ on $(\R^n, a)$, 
the {\it vanishing order} $\ord_a^\eta(h)$ of the function $h$ at the point $a$ for the vector-field $\eta$ is defined by 
$$
\ord_a^\eta(h) : = \inf \{ i \in \NN \cup \{ 0\} \mid (\eta^ih)(a) \not= 0 \}. 
$$
}
\enf

Then we characterize the cuspidal edge as an opening of fold map-germ: 

\bet
\label{Recognition-cuspidal-edge}
{\rm (Recognition of cuspidal edge)} 
For a frontal  $f : (\R^2, a) \to (\R^3, b)$ of corank $1$, 
the following conditions are equivalent to each other: 

{\rm (1)} $f$ is ${\mathcal A}$-equivalent to the cuspidal edge (CE). 

{\rm (1')} $f$ is ${\mathcal A}$-equivalent to the germ $(t_1, t_2) \mapsto (t_1, t_2^2, t_2^3)$. 

{\rm (2)} $f$ is a front and $\eta\lambda(a) \not= 0$. 

{\rm (3)} $\eta\lambda(a) \not= 0$ and $\ord_a^\eta(f_3) = 3$, 
for an adapted system of coordinates $(x_1, x_2, x_3)$ of $(\R^3, b)$. 
\ent

Theorem \ref{Recognition-cuspidal-edge} is generalized by 

\bet
\label{Recognition-cuspidal-edge2}
{\rm (Recognition of embedded cuspidal edge)} 
For a frontal  $f : (\R^2, a) \to (\R^m, b), 3 \leq m$ of corank $1$, 
the following conditions are equivalent to each other: 

{\rm (1)} $f$ is ${\mathcal A}$-equivalent to the cuspidal edge, i.e. the tangent surface to 
a curve of type $(1, 2, 3, \dots)$. 

{\rm (1')} $f$ is ${\mathcal A}$-equivalent to the germ $(t_1, t_2) \mapsto (t_1, t_2^2, t_2^3, 0, \dots, 0)$. 

{\rm (2)} $f$ is a front and $\eta\lambda(a) \not= 0$. 

{\rm (3)} $\eta\lambda(a) \not= 0$ and $\ord_a^\eta(f_i) = 3$ for some $i, 3 \leq i \leq m$, 
for an adapted system of coordinates $(x_1, x_2, x_3, \dots, x_m)$ of $(\R^m, b)$. 
\ent

The following is a recognition of the folded umbrella due to the theory of openings: 

\bet
\label{Recognition-FU}
{\rm (Recognition of folded umbrella (cuspidal cross cap))} 
Let $f : (\R^2, a) \to (\R^3, b)$ be a frontal of corank $1$. 
The following conditions are equivalent to each other: 

{\rm (1)} $f$ is ${\mathcal A}$-equivalent to the folded umbrella (FU), 
 i.e. the tangent surface to a curve of type $(1, 2, 4)$. 

{\rm (1')} $f$ is ${\mathcal A}$-equivalent to the germ $(t_1, t_2) \mapsto (t_1, \ t_2^2, \ t_1t_2^3)$. 

{\rm (2)} $\eta\lambda(a) \not= 0, (\eta^3 f_3)(a) = 0$ and $(d\lambda \wedge d(\eta^3 f_3))(a) \not= 0$. 
\ent

\ber
{\rm 
It is already known another kind of recognition of folded umbrella by \cite{FSUY}. 
%Fujimori-Saji-Umehara-Yamada (2008). 
}
\enr

As for cases of higher codimension, we have 

\bet
\label{OFU-characterization}
{\rm (Recognition of open folded umbrella (open cuspidal cross cap)).}
\\
Let $f : (\R^2, a) \to (\R^m, b), (m \geq 4)$ be a frontal of corank $1$. 
Then the following conditions are equivalent to each other: 

{\rm (1)} $f$ is ${\mathcal A}$-equivalent to the open folded umbrella, i.e. the 
tangent surface to a curve of type $(1, 3, 4, 5, \dots)$. 

{\rm (1')} $f$ is ${\mathcal A}$-equivalent to the germ $(t_1, t_2) \to 
(t_1, \ t_2^2, \ t_1t_2^3, \ t_2^5, \ 0, \dots, 0)$. 

{\rm (2)} $(\eta\lambda)(a) \not= 0$, $(\eta^3f_k)(a) = 0, (3 \leq k \leq m)$, 
and there exist $3 \leq i <  j \leq m$ and $A \in {\mbox{\rm GL}}(2, \R)$ 
such that, setting $(f_i, f_j)A = (f_3', f_4')$, $(d\lambda \wedge \eta^3f_3')(a) \not= 0, 
(d\lambda \wedge \eta^3f_ 4')(a) = 0, (\eta^5f_4')(a) \not= 0$. 
\ent

As for openings of Whitney's cusp mapping, we have

\bet
\label{Recognition-swallowtail}
{\rm (Recognition of swallowtail)} 
Let $f : (\R^2, a) \to (\R^3, b)$ be a frontal of corank $1$. 
Then the following conditions are equivalent to each other: 

{\rm (1)} $f$ is ${\mathcal A}$-equivalent to the swallowtail (SW), i.e. 
the tangent surface to a curve of type $(2, 3, 4)$. 

{\rm (1')} $f$ is ${\mathcal A}$-equivalent to the germ 
$(t_1, t_2) \mapsto (t_1, \ t_2^3 + t_1t_2, \ \frac{3}{4}t_2^4 + \frac{1}{2}t_1t_2^2)$. 

{\rm (2)} $f$ is a front, $(d\lambda)(a) \not= 0$ and $\ord_a^\eta(\lambda) = 2$. 

{\rm (3)} $\lambda$ is ${\mathcal K}$-equivalent to the germ $(t_1, t_2) \mapsto t_1$ at $0$, 
$\ord_a^\eta(\lambda) = 2$ and $\ord_a^\eta(f_3) = 4$, for an adapted system of coordinates $(x_1, x_2, x_3)$. 
\ent

As for cases of higher codimension, we have

\bet
\label{open-swallowtail}
{\rm  (Recognition of open swallowtail)}
Let $f : (\R^2, a) \to (\R^m, b)$ be a frontal of corank $1$ with $m \geq 4$. 
Then the following conditions are equivalent to each other: 

{\rm (1)} $f$ is ${\mathcal A}$-equivalent to the open swallowtail, i.e. 
the tangent surface to a curve of type $(2, 3, 4, 5, \dots)$. 

{\rm (1')} $f$ is ${\mathcal A}$-equivalent to the germ 
$
{\textstyle 
(t_1, t_2) \mapsto (t_1, \ t_2^3 + t_1t_2, \ \frac{3}{4}t_2^4 + \frac{1}{2}t_1t_2^2, \ \frac{3}{5}t_2^5 + \frac{1}{3}t_1t_2^3, 
\ 0, \dots). 
}
$

{\rm (2)} The Jacobian $\lambda$ is ${\mathcal K}$-equivalent to the germ $(t_1, t_2) \mapsto t_1$ at the origin, 
$\ord_a^\eta(\lambda) = 2$, $(\eta^3f_i)(a) = 0, (3 \leq k \leq m)$, 
and there exist $3 \leq i <  j \leq m$ and $A \in {\mbox{\rm GL}}(2, \R)$ 
such that, setting $(f_i, f_j)A = (f_3', f_4')$, 
$\ord_a^\eta(f_3') = 4$, $\ord_a^\eta(f_4') = 5$. 
\ent

\ber
{\rm 
Though we treat the open swallowtail as the singularity appeared in tangent surfaces, 
first it appeared as a singularity of Lagrangian varieties and geometric solutions of differential systems 
(\cite{Arnold}\cite{Givental}). 
The open swallowtail and open folded umbrella appear also in the context of 
frontal-symplectic versality (Example 12.3 of \cite{IJ}). 
}
\enr

As for openings of bec {\` a} bec mapping, we have

\bet
\label{Mond-singularities}
{\rm (Recognition of Mond singularity (cuspidal beaks),  (1)(2) \cite{IST}). }
Let $f : (\R^2, a) \to (\R^3, b)$ be a frontal of corank $1$. 
Then the following conditions are equivalent to each other: 

{\rm (1)} $f$ is ${\mathcal A}$-equivalent to Mond singularity (cuspidal beaks), 
i.e. 
the tangent surface to a curve of type $(1, 3, 4)$. 

{\rm (1')} $f$ is ${\mathcal A}$-equivalent to the germ 
$
{\textstyle 
(t_1, t_2) \mapsto (t_1, \ t_2^3 + t_1t_2^2, \ \frac{3}{4}t_2^4 + \frac{2}{3}t_1t_2^3). 
}
$

{\rm (2)} $f$ is a front, $\lambda$ is ${\mathcal K}$-equivalent $t_1t_2$ at the origin, 
and $\ord_a^\eta(\lambda) = 2$. 

{\rm (3)} $\lambda$ is ${\mathcal K}$-equivalent $t_1t_2$ at the origin, 
$\ord_a^\eta(\lambda) = 2$ and $\ord_a^\eta(f_3) = 4$. 
\ent

Moreover we have:

\bet
\label{open-Mond-singularities}
{\rm (Recognition of open Mond singularities (open cuspidal beaks)).}
Let $f : (\R^2, a) \to (\R^m, b)$ be a frontal of corank $1$ with $m \geq 4$. 
Then the following conditions are equivalent to each other: 

{\rm (1)} $f$ is ${\mathcal A}$-equivalent to the open Mond singularity, 
i.e. 
the tangent surface to a curve of type $(1, 3, 4, 5, \dots)$. 

{\rm (1')} 
$f$ is ${\mathcal A}$-equivalent to 
the germ 
$
{\textstyle 
(t_1, t_2) \mapsto (t_1, \ t_2^3 + t_1t_2^2, \ \frac{3}{4}t_2^4 + \frac{2}{3}t_1t_2^3, \ 
\frac{3}{5}t_2^5 + \frac{1}{2}t_1t_2^4, \ \dots). 
}
$

{\rm (2)} $\lambda$ is ${\mathcal K}$-equivalent to $(t_1, t_2) \mapsto t_1t_2$ at the origin, 
$\ord_a^\eta(\lambda) = 2$, 
$(\eta^3f_i)(a) = 0, (3 \leq k \leq m)$, 
and there exist $3 \leq i \not=  j \leq m$ and $A \in {\mbox{\rm GL}}(2, \R)$ 
such that, setting $(f_i, f_j)A = (f_3', f_4')$, 
$\ord_a^\eta(f_3') = 4$, $\ord_a^\eta(f_4') = 5$. 
\ent

To conclude this section, we give the result on recognition of Shcherbak singularity: 

\bet
\label{Shcherbak-singularity}
{\rm (Recognition of Shcherbak singularity)}
Let $f : (\R^2, a) \to (\R^3, b)$ be a frontal of corank $1$. 
Then the following conditions are equivalent to each other: 

{\rm (1)} $f$ is ${\mathcal A}$-equivalent to Shcherbak singularity, i.e. 
the tangent surface to a curve of type $(1, 3, 5)$. 

{\rm (1')} $f$ is ${\mathcal A}$-equivalent to the germ 
$(t_1, t_2) \mapsto (t_1, \ t_2^3 + t_1t_2^2, \ \frac{3}{5}t_2^5 + \frac{1}{2}t_1t_2^4)$ at the origin. 

{\rm (2)} $\lambda$ is ${\mathcal K}$-equivalent to the germ 
$(t_1, t_2) \mapsto t_1t_2$ at the origin, 
$\ord_a^\eta(\lambda) = 2$,  $\ord_c^\eta(f_3) \geq 4$ for any point 
$c$ on a component of the singular locus $S(f)$, and $\ord_a^\eta(f_3) = 5$. 
\ent

Note that Shcherbak singularity necessarily has $(2, 5)$ cuspidal-edge along one component of 
the singular locus, while it has ordinary $(2, 3)$ cuspidal edge along another component.

\section{Frontals and openings}
\label{Frontals and openings}

%\section{Jacobi modules}
%\label{Jacobi-kernel}

%To show a sequence of results on recognitions of frontal singularities, 
%we introduce a useful new notion 
%of map-germs (see also \cite{Ishikawa16}). 

To understand the frontal singularities and to prove the results in the previous section, 
we introduce the notion of openings and make clear its relation to frontal singularities (see also \cite{Ishikawa16}). 

Let $f : (\R^n, a) \to (\R^m, b)$ be a frontal (resp. a proper frontal) 
and $\widetilde{f} : (\R^n, a) \to {\mbox{\rm Gr}}(n, T\R^m)$ any Legendre lift of $f$. 
Let $(x_1, \dots, x_n, x_{n+1}, \dots, x_m)$ be an adapted system of coordinates 
to $\widetilde{f}$ (resp. to $f$) (Definition \ref{adapted-coordinates}). 
Then, setting $f_i = x_i\circ f, 1 \leq i \leq m$, we have 
$$
df_i = h_{i1}df_1 + h_{i2}df_2 + \cdots + h_{in}df_n,  \ (n+1 \leq i \leq m)
$$
for some $h_{ij} \in {\mathcal E}_a, h_{ij}(a) = 0$, 
$n+1 \leq i \leq m, 1 \leq j \leq n$. 

\bef
\label{Jacobi-module}
{\rm 
In general, for a map-germ $f = (f_1, \dots, f_m) : (\R^n, a) \to (\R^m, b)$, we define 
the ${\mathcal E}_a$-submodule 
$$
{\mathcal J}_f := \sum_{j=1}^m {\mathcal E}_{a} df_j = {\mathcal E}_a d(f^*{\mathcal E}_{b})
$$
of the ${\mathcal E}_a$-module of differential $1$-forms $\Omega_a^1$ on $(\R^n, a)$. 
We would like to call ${\mathcal J}_f$ the {\it Jacobi module} of $f$. 
}
\enf

Note that ${\mathcal J}_f$ is determined by the Jacobi matrix $J(f)$ of $f$. 
Returning to our original situation, we define the following key notion: 

\bef
{\rm 
We call a map-germ $f : (\R^n, a) \to (\R^m, b)$ an {\it opening} of a map-germ 
$g : (\R^n, a) \to (\R^n, g(a))$ if $f$ of form 
$$
(g_1, \dots, g_n, f_{n+1}, \dots, f_m)
$$ 
with 
$df_j \in {\mathcal J}_g, (n+1 \leq j \leq m)$ via a system of local coordinates of $(\R^m, b)$. 
}
\enf

Then we observe the following: 

\bel
\label{opening-frontal}
Any frontal $f : (\R^n, a) \to (\R^m, b)$ is an opening of 
$g := (f_1, \dots, f_n) : (\R^n, a) \to (\R^n, g(a))$ via adapted coordinates to a
Legendre lift of $f$. 
Conversely, any opening of a map-germ $g : (\R^n, a) \to (\R^n, g(a))$ is a frontal. 
An opening of $g$ is a proper frontal if and only if $g$ is proper, i.e. $S(g)$ is nowhere dense. 
\enl

\Proof
The first half is clear. To see the second half, 
let $f = (g_1, \dots, g_n, f_{n+1}, \dots, f_m)$ be an opening of $g$. 
Then 
$$
df_i = h_{i1}df_1 + h_{i2}df_2 + \cdots + h_{in}df_n,  \ (n+1 \leq i \leq m)
$$
for some $h_{ij} \in {\mathcal E}_a, n+1 \leq i \leq m, 1 \leq j \leq n$. 
Then a Legendre lift $\widetilde{f} : (\R^n, a) \to \Gr(n, T\R^m)$ is given, 
via Grassmannian coordinates of the fiber, by 
$$
t \mapsto (f(t), \ 
\left(
\begin{array}{c}
E_n \\
H(t)
\end{array}
\right)), 
$$
where $E_n$ is the $n\times n$ unit matrix and $H(t)$ is given by the $(m-n)\times n$-matrix $(h_{ij}(t))$. 
Therefore $f$ is a frontal. 
Note that an adapted system of coordinates for $f$ is given by 
$(x_1, \dots, x_n, \widetilde{x}_{n+1}, \dots, \widetilde{x}_{m})$ with 
$\widetilde{x}_{i} = x_{i} - \sum_{j=n+1}^m h_{ij}(a)x_j\  (n+1 \leq i \leq m)$. 

The last statement follows clearly. 
\QED

\

Here we recall one of key notion for our approach to the recognition problem of 
frontal singularities. 

\bef
\label{versal-opening}
{\rm 
(\cite{Ishikawa12}) 
An opening $f : (\R^n, a) \to (\R^m, b), f = (g; f_{n+1}, \dots, f_{m})$ of a map-germ 
$g : (\R^n, a) \to (\R^n, g(a))$ is called a {\it versal opening} if, 
for any $h \in {\mathcal E}_a$ with $dh \in {\mathcal J}_g$, 
there exist $k_0, k_1, \dots, k_{m-n} \in {\mathcal E}_{\R^n, g(a)}$ 
such that 
$$
h = g^*(k_0) + g^*(k_1)f_{n+1} + \cdots + g^*(k_{m-n})f_m. 
$$
}
\enf

We will use the following result which is proved in Proposition 6.9 of \cite{Ishikawa12}. 

\bet
\label{uniqueness-versal-opening}
Any two versal openings $f, f' : (\R^n, a) \to (\R^m, b)$ (having the same target dimension) 
of a map-germ $g$ are ${\mathcal A}$-equivalent to each other. 
\ent

Recall, for a map-germ $f : (\R^n, a) \to (\R^m, b)$, we have defined 
${\mathcal J}_f = {\mathcal E}_{a} d(f^*{\mathcal E}_{b})$ (Definition \ref{Jacobi-module}). 
%
%First we show

\bel
\label{right-left-J}
{\rm (1)} Let $f : (\R^n, a) \to (\R^m, b)$, $f' : (\R^n, a) \to (\R^{m}, b')$ be map-germs. 
If $f$ and $f'$ are ${\mathcal L}$-equivalent, i.e. if there exists a diffeomorphism-germ 
$\tau : (\R^m, b) \to (\R^m, b')$ such that $f' = \tau\circ f$, then 
${\mathcal J}_f = {\mathcal J}_{f'}$. 

{\rm (2)} Let $f : (\R^n, a) \to (\R^m, b)$, $f' : (\R^n, a') \to (\R^{m}, b)$ be map-germs. 
If $f$ and $f'$ are ${\mathcal R}$-equivalent, i.e. if there exists a diffeomorphism-germ 
$\sigma : (\R^n, a) \to (\R^n, a')$ such that $f' = f\circ\sigma$, then 
$\sigma^*({\mathcal J}_f) = {\mathcal J}_{f'}$. 
\enl

\Proof
(1) Since $f^*{\mathcal E}_b = {f'}^*{\mathcal E}_{b'}$, we have 
${\mathcal J}_f = {\mathcal E}_a d(f^*{\mathcal E}_b) =  {\mathcal E}_a d({f'}^*{\mathcal E}_{b'}) 
= {\mathcal J}_{f'}$. 

(2) Since ${f'}^*{\mathcal E}_b = \sigma^*(f^*{\mathcal E}_b)$, we have 
$$
{\mathcal J}_{f'} = {\mathcal E}_{a'} d({f'}^*{\mathcal E}_b) 
= {\mathcal E}_{a'} d(\sigma^*(f^*{\mathcal E}_b)) = \sigma^*{\mathcal E}_{a} \sigma^*d(f^*{\mathcal E}_b)
= \sigma^*({\mathcal E}_{a} d(f^*{\mathcal E}_b)) = \sigma^*({\mathcal J}_f). 
$$
\QED

The equality of Jacobi modules ${\mathcal J}_f$ has a simple meaning:

\bel
\label{J-equivalence-lemma}
Let $f : (\R^n, a) \to (\R^m, b)$, $f' : (\R^n, a) \to (\R^{m'}, b')$ be map-germs. 

Then the following conditions {\rm (i), (ii)} are equivalent: 

{\rm (i)} The Jacobi module ${\mathcal J}_f = {\mathcal J}_{f'}$. 

{\rm (ii)} There exist an $m'\times m$-matrix $P$ and an $m\times m'$-matrix $Q$ 
with entries in ${\mathcal E}_a$ such that 
the Jacobi matrix $J(f') = P J(f)$ and $J(f) = Q J(f')$. 

In particular, {\rm (i)} implies that the Jacobi ideal $J_f = J_{f'}$. 

Moreover, if the target dimension $m = m'$, then the following condition (iii) is equivalent to (i). 

{\rm (iii)} There exists an invertible $m\times m$-matrix $R$ with entries in ${\mathcal E}_a$ such that 
$J(f') = R J(f)$. 
\enl

To show Lemma \ref{J-equivalence-lemma}, we recall the following fact in linear algebra. 

\bel
\label{Mather-contact}
{\rm (cf. \cite{Mather})} 
Let $A, B$ be $m\times m$-matrices with entries in $\R$. Then there exists an $m\times m$-matrices 
$C$ with entries in $\R$ such that $C(E_m - BA) + A$ is invertible. 
\enl

\noindent
{\it 
Proof of Lemma \ref{J-equivalence-lemma}:} 

The inclusion ${\mathcal J}_{f'} \subseteq {\mathcal J}_{f}$ is equivalent 
to that there exist $p_{ij} \in {\mathcal E}_a$ such that 
$df_i' = \sum_{j=1^m} p_{ij} df_j, (1 \leq i \leq m)$, namely that 
$J(f') = P J(f)$ by setting $P = (p_{ij})$. 
Similarly, the inclusion 
${\mathcal J}_{f} \subseteq {\mathcal J}_{f'}$ 
is equivalent to that there exist $q_{ij} \in {\mathcal E}_a$ such that 
$df_i = \sum_{j=1^m} q_{ij} df_j', (1 \leq i \leq m)$, namely that 
$J(f) = Q J(f')$ by setting $Q = (q_{ij})$. 
Therefore the equivalence between (i) and (ii) is clear. 

Suppose $m = m'$. 
By Lemma \ref{Mather-contact}, there exists an $m\times m$-matrix $C$ with entries in $\R$ 
such that $C(E_m - Q(a)P(a)) +P(a)$ is invertible. 
Then $R := C(E_m - QP) +P$ is an invertible $m\times m$-matrix with entries in 
in ${\mathcal E}_a$. 
Then we have $(E_m - QP)J(f) = J(f) - Q J(f') = O$ and therefore 
$R J(f) = C(E_m - QP)J(f) + P J(f) = J(f')$. 
\QED

\

\ber
\label{ramification} 
{\rm 
Related to Jacobi modules, we define the ramification module ${\mathcal R}_f \subseteq {\mathcal E}_a$ 
for a map-germ $f : (\R^n, a) \to (\R^m, b)$ by
$$
{\mathcal R}_f := \{ h \in {\mathcal E}_a \mid dh \in {\mathcal J}_f\}, 
$$
using the Jacobi module ${\mathcal J}_f$. Then ${\mathcal R}_f = {\mathcal R}_{f'}$ 
if and only if ${\mathcal J}_f = {\mathcal J}_{f'}$. See, for details, the series of papers 
\cite{Ishikawa95}\cite{Ishikawa96}\cite{Ishikawa12}\cite{Ishikawa13}\cite{Ishikawa14}\cite{Ishikawa16}. 
}
\enr

\bel
Let $f : (\R^n, a) \to (\R^m, b)$, $f' : (\R^n, a') \to (\R^{m'}, b')$ be map-germs. 
If 
${\mathcal J}_f = {\mathcal J}_{f'}$, then 
$$
J_f = J_{f'}, \quad 
{\mathcal N}_f = {\mathcal N}_{f'}. 
$$
\enl

\Proof
The equality $J_f = J_{f'}$ follows from 
Lemma \ref{J-equivalence-lemma}. 
For any $\eta \in {\mathcal V}_a$,  the condition $\eta \in {\mathcal N}_f$ is 
equivalent to that $\omega(\eta) \in J_f = J_{f'}$ 
for any ${\mathcal J}_f = {\mathcal J}_{f'}$, which is equivalent to that $\eta \in {\mathcal N}_{f'}$. 
Therefore we have ${\mathcal N}_f = {\mathcal N}_{f'}$. 
\QED

\

\bel
\label{corank-one-converse}
Let $f, f' : (\R^n, a) \to (\R^m, b)$ be a proper frontal of corank $1$. Then the conditions 
$$
\lambda_f\cdot{\mathcal E}_a = \lambda_{f'}\cdot{\mathcal E}_a, \quad 
{\mathcal N}_f = {\mathcal N}_{f'}, 
$$
imply that ${\mathcal J}_f = {\mathcal J}_{f'}$. 
\enl

\Proof
By the assumption we may take $\lambda_f = \lambda_{f'}$ and $\eta_f = \eta_{f'}$. 
and $\eta_f = \pa/\pa t_n$ for a system of coordinates 
$t_1, \dots, t_{n-1}, t_n$ of $(\R^n, a)$. 
Note that, by the assumption, the zero-locus of $\lambda_f$ is nowhere dense. 
Then 
$f_*(\pa/\pa t_1), \dots, f_*(\pa/\pa t_{n-1}), (1/\lambda_f)f_*(\pa/\pa t_n)$ 
are linearly independent at $a$ as elements of ${\mathcal E}_a^m$. 
Take additional $\xi_{n+1}, \dots, \xi_m$ to complete a basis of ${\mathcal E}_a^m$. 
Moreover by the assumption 
$$
{f'}_*(\pa/\pa t_1), \dots, {f'}_*(\pa/\pa t_{n-1}), (1/\lambda_f){f'}_*(\pa/\pa t_n)
$$ 
are linearly independent at $a$ as elements of ${\mathcal E}_a^m$. 
Take additional $\xi'_{n+1}, \dots, \xi'_m$ to complete a basis of ${\mathcal E}_a^m$. 
Then define $R : (\R^n, a) \to \GL(m, \R)$ by 
$$
Rf_*(\pa/\pa t_i) = {f'}_*(\pa/\pa t_i), 1 \leq i \leq n-1, \ \ 
R(1/\lambda_f)f_*(\pa/\pa t_n) = (1/\lambda_f){f'}_*(\pa/\pa t_n), \ \ 
R\xi_j = \xi_j', n+1 \leq j \leq m. 
$$
Then $R f_*(\pa/\pa t_n) = {f'}_*(\pa/\pa t_n)$ and we have $R J(f) = J(f')$. 
By Lemma \ref{J-equivalence-lemma}, we have ${\mathcal J}_f = {\mathcal J}_{f'}$. 
\QED

\

We utilize the following in the next section: 

\bel
\label{opening-Jacobian-kernel}
Let $f : (\R^n, a) \to (\R^m, b)$ be an opening of $g : (\R^n, a) \to (\R^n, g(a))$ 
with respect to an adapted system of coordinates $(x_1, \dots, x_n, x_{n+1}, \dots, x_m)$. 
Then $f$ and $g$ are frontals and ${\mathcal J}_f = {\mathcal J}_g$. They have common Jacobian, 
same corank, and ${\mathcal N}_f = {\mathcal N}_g$. If they are of corank $1$, then 
they have common kernel field. 
\enl

\Proof
By Lemma \ref{opening-frontal}, we have 
${\mathcal J}_f = {\mathcal J}_g$. 
Then $J_f = J_g$, therefore $\lambda_f = \lambda_g$. 
Moreover, by Lemma \ref{J-equivalence-lemma}, $\Ker(T_a f) = \Ker(T_a g) \subseteq T_a\R^n$. 
Therefore $f$ and $g$ have the same corank. Furthermore, for any $\eta \in {\mathcal V}_a$, the condition 
that $df_i(\eta) \in J_f, 1 \leq i \leq m$ is equivalent to that $dg_i(\eta) \in J_f = J_g, 1 \leq i \leq n$. 
Hence ${\mathcal N}_f = {\mathcal N}_g$. 
\QED

\bef
\label{J-equivalence}
{\rm 
Let $f : (\R^n, a) \to (\R^m, b)$ and $f' : (\R^n, a') \to (\R^{m'}, b')$ be map-germs. 
Then $f$ and $f'$ are called {\it ${\mathcal J}$-equivalent} if 
there exists a diffeomorphism-germ $\sigma : (\R^n, a) \to (\R^n, a')$ such that 
${\mathcal J}_{f'\circ\sigma} = {\mathcal J}_f$. 
Note that $m$ and $m'$ can be different. 
}
\enf

By Lemma \ref{right-left-J} and Lemma \ref{corank-one-converse}, we have 

\bec
Let $f : (\R^n, a) \to (\R^m, b)$ and $f' : (\R^n, a') \to (\R^{m'}, b')$ be map-germs. 
If $f$ and $f'$ are ${\mathcal A}$-equivalent, then $f$ and $f'$ are ${\mathcal J}$-equivalent. 
\enc

\bec
Let $f, f'$ be proper frontals. If
$f$ and $f'$ are ${\mathcal J}$-equivalent, then 
$(\lambda_f\cdot {\mathcal E}_{a}, \, {\mathcal N}_f)$ is 
transformed to 
$(\lambda_{f'}\cdot{\mathcal E}_{a'}, \, {\mathcal N}_{f'})$ by a diffeomorphism-germ 
$\sigma : (\R^n, a) \to (\R^n, a')$. In particular $\lambda_f$ and $\lambda_{f'}$ are 
${\mathcal K}$-equivalent. 

Moreover if $f$ is of corank $1$ and 
$(\lambda_f\cdot {\mathcal E}_{a}, \, {\mathcal N}_f)$ is 
transformed to 
$(\lambda_{f'}\cdot{\mathcal E}_{a'}, \, {\mathcal N}_{f'})$ by a diffeomorphism-germ 
$\sigma : (\R^n, a) \to (\R^n, a')$, then $f$ and $f'$ are ${\mathcal J}$-equivalent. 
\enc

On the vanishing order of function for a vector field introduced in Definition \ref{vanishing-order}, we have: 

\bel
\label{order-invariance}
If $\widetilde{h} = \rho h, \widetilde{\xi} = \nu \xi$ for some $\rho, \nu \in {\mathcal E}_a$ 
with $\rho(a) \not= 0, \xi(a) \not= 0$, then 
$\ord_a^{\widetilde{\xi}}(\widetilde{h}) = \ord_a^\xi(h)$. 
If $\overline{h} =  h\circ\sigma, \overline{\xi} = (T\sigma^{-1})\circ\xi\circ\sigma$ 
for some $\sigma : (\R^n, a') \to (\R^n, a)$ is a diffeomorphism-germ, then 
$\ord_{a'}^{\overline{\xi}}(\overline{h}) = \ord_a^\xi(h)$. 
\enl

By Lemma \ref{order-invariance} we have 

\bec
Let $f : (\R^n, a) \to (\R^m, b)$ be a proper frontal of corank $1$. Then 
$\ord_a^\eta(\lambda)$ is independent of the choices of the Jacobian $\lambda$ and the kernel 
field $\eta$ of $f$. If $f' : (\R^n, a') \to (\R^{m'}, b')$ is ${\mathcal J}$-equivalent to $f$, then 
$f'$ is a proper frontal of corank $1$ and $\ord_a^{\eta'}(\lambda')$ is equal to $\ord_a^\eta(\lambda)$, 
for any Jacobian $\lambda'$ and any kernel filed $\eta'$ of $f'$. 
\enc

\section{Proofs of recognition theorems} 
\label{Proofs-of-recognition-theorems} 

In this section we give proofs of Theorems \ref{Recognition-cuspidal-edge}, \ref{Recognition-cuspidal-edge2}, \ref{Recognition-FU}, 
\ref{OFU-characterization}, \ref{Recognition-swallowtail}, \ref{open-swallowtail}, \ref{Mond-singularities}, 
\ref{open-Mond-singularities}, and \ref{Shcherbak-singularity}. 

\

\noindent
{\it Proof of Theorem \ref{Recognition-cuspidal-edge}:} 
The equivalence of (1) and (1') is classically known (see \cite{Ishikawa95}). 
The equivalence of (1') and (2) is proved in \cite{KRSUY}. 

To study the condition, we set $g = (f_1, f_2)$. Then for the Jacobian $\lambda$ and the kernel field 
$\eta$ of $g$ we also have $\eta\lambda(a) \not= 0$ (see Lemma\ref{opening-Jacobian-kernel}). 
By Theorem \ref{Whitney, Saji} $g$ is ${\mathcal A}$-equivalent to fold. 
Then the condition (3) means that $f$ is a versal opening of the fold $g$. 
Since the cuspidal edge is characterized as the (mini)-versal opening of the fold mp-germ, we have 
the equivalence of (3) and (1) by Theorem \ref{uniqueness-versal-opening}. 
\QED

\

\noindent
{\it Proof of Theorem \ref{Recognition-cuspidal-edge2}:} 
The equivalence of (1) and (1') is proved in Theorem 7.1 of \cite{Ishikawa12}. 
The condition (3) means that $f$ is a versal opening of the fold $g$. 
Since the embedded cuspidal edge is characterized as the versal opening of the fold mp-germ, we have 
the equivalence of (3) and (1) by Theorem\ref{uniqueness-versal-opening}. 
On the other hand, under the condition $\eta\lambda(a) \not= 0$, 
the condition $\ord_a^\eta(f_i) = 3$ for some $i, 3 \leq i \leq m$ is equivalent to that the Legendre lift $\widetilde{f}$ 
is an immersion i.e. $f$ is a front. Therefore (3) and (2) are equivalent. 
\QED

\

\noindent
{\it Proof of Theorem \ref{Recognition-FU}.} 
The equivalence of (1) and (1') is due to Cleave (see \cite{Ishikawa12}). 

Suppose the condition (2) is satisfied. Then $f$ is ${\mathcal A}$-equivalent to the germ 
$g(t_1, t_2) = (t_1, t_2^2, f_3(t_1, t_2))$ at the origin with $\lambda = t_2, \eta = \pa/\pa t_2$, 
$(\eta^3 f_3)(0) = 0$ and $(d\lambda \wedge d(\eta^3 f_3))(0) \not= 0$. 
Since $df_3 \in {\mathcal J}_g$, in other word since $f_3 \in {\mathcal R}_g$ (Remark\ref{ramification}), 
there exist functions $A, B$ on $(\R^2, 0)$ such that 
$$
f_3(t_1, t_2) = A(t_1, t_2^2) + B(t_1, t_2^2) t_2^3. 
$$
Then the condition $(\eta^3 f_3)(0) = 0$ is equivalent to $B(0, 0) = 0$, 
and the condition $(d\lambda \wedge d(\eta^3 f_3))(0) \not= 0$ is 
equivalent to $\frac{\pa B}{\pa t_1}(0, 0) \not= 0$. 
Define diffeomorphism-germs 
$\sigma : (\R^2, 0) \to (\R^2, 0)$ by $\sigma(t_1, t_2) = (B(t_1, t_2^2), t_2)$ 
and $\tau : (\R^3, 0) \to (\R^3, 0)$ by 
$\tau(x_1, x_2, x_3) = (B(x_1, x_2), x_2, x_3 - A(x_1, x_2))$. 
Then $(t_1, t_2^2, t_1t_2^3)\circ\sigma = \tau\circ(t_1, t_2^2, f_3)$ holds. Therefore 
$f$ is ${\mathcal A}$-equivalent to folded umbrella. Hence we see that 
(2) implies (1). Conversely (1) implies (2) for some, so for any, adapted coordinates. 
\QED

\

\noindent
{\it Proof of Theorem \ref{OFU-characterization}:} 
The ${\mathcal A}$-determinacy of tangent maps to curves of type $(1, 2, 4, 5, \dots)$ 
is proved in Theorem 7.2 of \cite{Ishikawa12}. 
Let $\gamma : (\R, 0) \to (\R^m, 0)$ be the curve 
$t \mapsto (t, t^2, t^4, t^5, 0, \dots)$. Then the tangent map $\Tan(\gamma) : (\R^2, 0) 
\to (\R^m, 0)$ is given by 
$$
\Tan(\gamma)(t, u) = (t+u, \ t^2 + 2ut, \ t^4 + 4ut^3, \ t^5 + 5ut^4, \ 0, \ \dots). 
$$
Then it is easy to see that $\Tan(\gamma)$ is ${\mathcal A}$-equivalent to 
$(t_1, t_2) \to (t_1, \ t_2^2, \ t_1t_2^3, \ t_2^5, \ 0, \dots, 0)$. Hence we have 
the equivalence of (1) and (1'). 

Suppose $f$ satisfies (2). 
Then $f$ is an opening of $(f_1, f_2)$, which is a fold by Theorem\ref{Whitney, Saji}. 
Therefore $f$ is ${\mathcal A}$-equivalent to a frontal of form
$(t_1, t_2^2, f_3, f_4, \dots)$ for an adapted coordinates. The Jacobian is given by $\lambda = t_2$ 
and the kernel field is given by $\eta = \pa/\pa t_2$. 
We write $f_i = A_i(t_1, t_2^2) + B_i(t_1, t_2^2)t_2^3$ for some 
$A_i, B_i$ with $A_i(0, 0) = 0, B_i(0, 0) = 0, \ (3 \leq i \leq m)$. 
Then $f_i = \widetilde{A}_i(t_1, t_2^2)t_1t_2^3 + \widetilde{B}_i(t_1, t_2^2)t_2^5$. 
Then the condition (2) is equivalent to that, for some $i, j$ with $3 \leq i < j \leq m$,  
$$
\left(
\begin{array}{cc}
\widetilde{A}_i(0, 0) & \widetilde{B}_i(0, 0)
\\
\widetilde{A}_j(0, 0) & \widetilde{B}_j(0, 0)
\end{array}
\right) \in \GL(2, \R). 
$$
Then $f$ is ${\mathcal A}$-equivalent to $(t_1, \ t_2^2, \ t_1t_2^3, \ t_2^5, \ 0, \dots, 0)$. 
Therefore (2) implies (1'). The converse is clear. 
\QED

\

\noindent
{\it Proof of Theorem \ref{Recognition-swallowtail}:} 
The equivalence of (1) and (1') is proved in Theorem 1 of \cite{Ishikawa95}. 
The equivalence of (1') and (2) is proved in Proposition 1.3 of \cite{KRSUY}. 
The condition that $\lambda$ is ${\mathcal K}$-equivalent to $t_1$ and $\ord_a^\eta(\lambda) = 2$ 
is equivalent, by Theorem \ref{Whitney, Saji}, to that $f$ is an opening of Whitney's cusp 
$g(t_1, t_2) = (t_1, \ t_2^3 + t_1t_2)$. The Jacobian is given by $\lambda = 3t_2^2 + t_1$ and 
the kernel field is given by $\eta = \pa/\pa t_2$. 
Set $U_1 = \frac{3}{4}t_2^4 + \frac{1}{2}t_1t_2^2, U_2 = \frac{3}{5}t_2^5 + \frac{1}{3}t_1t_2^3$. 
Then it is known that the ramification module ${\mathcal R}_g$ is generated 
by $1, U_1, U_2$ over $g^*$ (see \cite{Ishikawa95}). 
Since $f_3 \in {\mathcal R}_g$ is the third component for an adapted system of coordinates, 
$f_3$ is written as 
$f_3 = A\circ g + (B\circ g)U_1 + (C\circ g)U_2$, 
for some functions $A, B, C$ with $A(0, 0) = 0, \frac{\pa A}{\pa x_1}(0, 0) = 0, \frac{\pa A}{\pa x_2}(0, 0) = 0$. 
By the condition $\ord_a^\eta(f_3) = 4$, we have $B(0, 0) \not= 0$. 
Then, by a change of adapted system of coordinates, 
We may suppose $f = (g, f_3)$ with $f_3 = U_1 + \Phi$, where 
$\Phi = (\widetilde{B}\circ g)U_1 + (D\circ g)U_2$ with $\widetilde{B}(0, 0) = 0$. Then 
we set the family $F_s = (g, U_1 + s\Phi$. By the same infinitesimal method used in \cite{Ishikawa95}, 
we can show that the family $F_s$ is trivialized by ${\mathcal A}$-equivalence. 
Hence $f = F_1$ is ${\mathcal A}$-equivalent to $F_0$, that is the normal form of (2). 
Therefore (3) implies (2). The converse is clear. 
\QED

\

\noindent
{\it Proof of Theorem \ref{open-swallowtail}:} 
The equivalence of (1) and (1') is proved in \cite{Ishikawa12}. 
The condition (2) implies, by Theorem \ref{Whitney, Saji}, that $f$ is an opening of Whitney's cusp. 
Using the same notations as in the proof of Theorem \ref{Recognition-swallowtail}, 
we write $f_k$ as 
$f_k = A_k\circ g + (B_k\circ g)U_1 + (C_k\circ g)U_2$, 
for some functions $A_k, B_k, C_k$ with $A_k(0, 0) = 0, \frac{\pa A_k}{\pa x_1}(0, 0) = 0, 
\frac{\pa A_k}{\pa x_2}(0, 0) = 0$. 
Then by the condition (2), we see that $f$ is a versal opening (Definition \ref{versal-opening}) of $g$. 
On the other hand the map-germ of (1') is a versal opening of $g$ (\cite{Ishikawa12}). 
By Theorem \ref{uniqueness-versal-opening}, we see that 
(2) implies (1'). The converse implication (1') to (2) is clear. 
\QED

\

\noindent
{\it Proof of Theorem \ref{Mond-singularities}:} 
The outline of the proof is similar to that of Theorem \ref{Recognition-swallowtail}. 
The equivalence of (1) and (1') is proved in Theorem 1 of \cite{Ishikawa95}. 
The equivalence of (1') and (2) is proved in \cite{IST}. 
The condition that $\lambda$ is ${\mathcal K}$-equivalent to $t_1t_2$ and $\ord_a^\eta(\lambda) = 2$ 
is equivalent, by Theorem \ref{Whitney, Saji}, to that $f$ is an opening of bec {\` a} bec 
$g(t_1, t_2) = (t_1, \ t_2^3 + t_1t_2^2)$. The Jacobian is given by $\lambda = 3t_2^2 + 2t_1t_2$ and 
the kernel field is given by $\eta = \pa/\pa t_2$. 
Set $U_1 = \frac{3}{4}t_2^4 + \frac{2}{3}t_1t_2^3, U_2 = \frac{3}{5}t_2^5 + \frac{1}{2}t_1t_2^4$. 
Then it is known that the ramification module ${\mathcal R}_g$ is generated 
by $1, U_1, U_2$ over $g^*$ (see \cite{Ishikawa95}). 
Since $f_3 \in {\mathcal R}_g$ is the third component for an adapted system of coordinates, 
$f_3$ is written as 
$f_3 = A\circ g + (B\circ g)U_1 + (C\circ g)U_2$, 
for some functions $A, B, C$ with $A(0, 0) = 0, \frac{\pa A}{\pa x_1}(0, 0) = 0, \frac{\pa A}{\pa x_2}(0, 0) = 0$. 
By the condition $\ord_a^\eta(f_3) = 4$, we have $B(0, 0) \not= 0$. 
Then, by a change of adapted system of coordinates, 
We may suppose $f = (g, f_3)$ with $f_3 = U_1 + \Phi$, where 
$\Phi = (\widetilde{B}\circ g)U_1 + (C\circ g)U_2$ with $\widetilde{B}(0, 0) = 0$. Then, 
by the infinitesimal method used in \cite{Ishikawa95}, 
the family $F_s = (g, U_1 + s\Phi$ is trivialized by ${\mathcal A}$-equivalence. 
Hence $f = F_1$ is ${\mathcal A}$-equivalent to $F_0$, that is the normal form of (2). 
Therefore (3) implies (2). The converse is clear. 
\QED

\

\noindent
{\it Proof of Theorem \ref{open-Mond-singularities}:} 
Open Mond singularities are characterized as versal openings of bec {\` a} bec (\cite{Ishikawa12}). 
Then Theorem\ref{open-Mond-singularities} is proved similarly as the proof of 
Theorem\ref{open-swallowtail}. 
\QED

\

\noindent
{\it Proof of Theorem \ref{Shcherbak-singularity}:} 
The equivalence of (1) and (1') is proved in \cite{Ishikawa95}. The condition (2) implies that 
$f$ is an opening of bec {\` a} bec. 
Using the same notations in the proof of Theorem \ref{Mond-singularities}, 
we write $f_3$ as 
$f_3 = A\circ g + (B\circ g)U_1 + (C\circ g)U_2$, 
for some functions $A, B, C$ with $A(0, 0) = 0, \frac{\pa A}{\pa x_1}(0, 0) = 0, \frac{\pa A}{\pa x_2}(0, 0) = 0$. 
By the condition $\ord_a^\eta(f_3) = 5$, we have $B(0, 0) = 0$ and $C(0, 0) \not= 0$. 
Moreover, by the assumption, we may assume that $\ord^\eta_{(t_1, 0)}f_3 \geq 4$ along the 
component $\{ t_2 = 0\}$ of $S(f)$ and then $B(x_1, 0) = 0$. 
Then, by a change of adapted system of coordinates, 
we may suppose $f = (g, f_3)$ with $f_3 = U_2 + \Phi$, where 
$\Phi = (B\circ g)U_1 + (\widetilde{C}\circ g)U_2$ with $B(x_1, 0) = 0, \widetilde{C}(0, 0) = 0$. Then 
by the same infinitesimal method used in \cite{Ishikawa95}, the family $F_s = (g, U_2 + s\Phi)$ turns to be 
trivial under ${\mathcal A}$-equivalence. 
Hence $f = F_1$ is ${\mathcal A}$-equivalent to $F_0$, that is the normal form of (1'). 
Therefore (2) implies (1'). The converse is clear. 
\QED

\section{An application to $3$-dimensional Lorentzian geometry, and other topics}
\label{application}

We announce the following result without explanations of notions. 
The details will be given in \cite{IMT18}. 

\bet
\label{Lorentzian-}
{\rm 
(\cite{CI}, \cite{IMT11}\cite{IMT18}) }
Any null frontal surface in a Lorentzian $3$-manifold turns to be a null tangent surface of a (directed) null curve, and 
any generic null frontal surface has only singularities, along the null curve, 

{\rm (I)} cuspidal edge {\rm (CE)}, 
{\rm (II)} swallowtail {\rm (SW)}, or {\rm (III)} Shcherbak singularity {\rm (SB)}. 

Moreover the corresponding dual frontal in the space of null-geodesics has 
{\rm (I)} cuspidal edge {\rm (CE)}, {\rm (II)} Mond singularity {\rm (MD)}, or {\rm (III)} generic folded pleat {\rm (GFP)}. 

The same classification result holds not only for any Lorentzian metric but also for arbitrary 
non-degenerate (strictly convex) cone structure in any $3$-manifold. 
\ent

To show Theorem \ref{Lorentzian-}, we face the recognition problem on cuspidal edge, swallowtail, 
Scherbak singularity, Mond singularity, and \lq\lq generic folded pleat". In fact we will use 
the recognition theorems introduced in the previous section and the following 
result on openings of Whitney's cusp. 
The following recognition result is proved 
by the same method of the above proof of Theorem\ref{Recognition-swallowtail}. 
The details will be given in \cite{IMT18}. 

\bet
\label{folded-pleat-recognition}
{\rm  (Recognition of folded pleat) }
Let $f : (\R^2, a) \to (\R^3, b)$ be a frontal of corank $1$. 
Then the following conditions are equivalent to each other: 

{\rm (1)} $f$ is ${\mathcal A}$-equivalent to a folded pleat i.e. 
the singularity of tangent surface of a curve of type $(2, 3, 5)$. 

{\rm (1')} $f$ is ${\mathcal A}$-equivalent to 
the germ $(t_1, \ t_2^3 + t_1t_2, \  \frac{3}{5}t_2^5 + \frac{1}{2}t_1t_2^3 + c(\frac{1}{2}t_2^6 + \frac{3}{4}t_1t_2^4))$
at the origin for some $c \in \R$. 

{\rm (2)} 
$\lambda$ is ${\mathcal K}$-equivalent to the germ $(t_1, t_2) \mapsto t_1$ at the origin, 
$\ord_a^\eta(\lambda)(a) = 2$, $f$ has an injective representative, and $\ord^\eta_p(f_3) = 5. $
\ent

Note that a folded pleat singularity necessarily has an injective representative.

%\begin{center}
%     \includegraphics[width=2.5truecm,height=2truecm, clip, 
%     bb=180 588 415 792]{folded-pleats18.pdf}
%     \hspace{1.5truecm} 
%     \includegraphics[width=2.5truecm,height=2truecm, clip, 
%     bb=77 706 236 816]{cuspidal-swallowtail.pdf}
%     \hspace{1.5truecm} 
%     \includegraphics[width=2.5truecm,height=2truecm, clip, 
%     bb=151 317 388 545]{cuspidal-lips.pdf}
%     \\
%     folded pleat  \hspace{1.5truecm} cuspidal swallowtail \hspace{1.8truecm} cuspidal lips
%\end{center}

\ber
{\rm 
Recall that the diffeomorphism classes (CE), (SW), (SB) and (MD) are exactly characterized as those of tangent surfaces 
in Euclidean space $\R^3$ of curves of 
type $(1, 2, 3)$, $(2, 3, 4)$, $(1, 3, 5)$, $(1, 3, 4)$ respectively. 
A map-germ $(\R^2, a) \to (\R^3, b)$ is called a {\it folded pleat} (FP) if it is diffeomorphic to the tangent surface 
of a curve of type $(2, 3, 5)$ in $\R^3$. 
The diffeomorphism classes of folded pleats fall into {\it two} classes,  
the generic folded pleat and the non-generic folded pleat. 
In the list of Theorem \ref{Lorentzian-}, it is claimed that only the generic folded pleat (GFP) appear. 
Therefore Theorem\ref{folded-pleat-recognition} 
do not solve the recognition of a singularity but a class of singularities, which 
consists of two singularities. Note that the parameter $c$ in (1') of Theorem\ref{folded-pleat-recognition} 
is not a moduli, but provides just two ${\mathcal A}$-equivalence classes. 
To recognize the generic folded pleat, it is necessary 
an additional argument to distinguish generic and non-generic folded pleats. 
}
\enr

\

In this occasion we introduce and prove the following two theorems of recognition:

\bet
{\rm (Recognition of cuspidal swallowtail)}\ 
Let $(\R^2, a) \to (\R^3, b)$ be a frontal of corank $1$. 
Then the following conditions are equivalent to each other: 

{\rm (1)} $f$ is ${\mathcal A}$-equivalent to the cuspidal swallowtail 
i.e. the singularity of tangent surface of curves of type $(3, 4, 5)$. 

{\rm (1')} $f$ is ${\mathcal A}$-equivalent to 
the germ 
$(t_1, t_2) \mapsto (t_1, \ t_2^4 + t_1t_2, \ \frac{4}{5}t_2^5 + \frac{1}{2}t_1t_2^2)$ at the origin. 

{\rm (2)} $\lambda$ is ${\mathcal K}$-equivalent to the germ 
$(t_1, t_2) \mapsto t_1$ at the origin, 
$\ord_a^\eta(\lambda) = 3$ and $\ord_a^\eta(f_3) = 5$. 
\ent

\Proof
In \cite{Ishikawa12} it is proved that the condition (1) is equivalent to that $f$ is ${\mathcal A}$-equivalent 
to the germ $(t, u) \mapsto (t^3 + 3u, t^4 + 4ut, t^5 + 5ut^2)$, which is ${\mathcal A}$-equivalent 
to the normal form of (1'). Therefore (1) and (1') are equivalent. 
In \cite{Saji}, the map-germ which is ${\mathcal A}$-equivalent to the germ 
$g : (t_1, t_2) \mapsto (t_1, \ t_2^4 + t_1t_2)$ at the origin is called a {\it swallowtail} and 
it is shown that a map-germ $g : (\R^2, a) \to (\R^2, g(a))$ is a swallowtail if and only if 
$\lambda$ is ${\mathcal K}$-equivalent to the germ 
$(t_1, t_2) \mapsto t_1$ at the origin and 
$\ord_a^\eta(\lambda) = 3$. 
Suppose $f$ satisfies (2). Then $f$ is an opening of swallowtail. 
Then $f$ is ${\mathcal A}$-equivalent 
to a frontal of form $f = (g, f_3)$. We have the Jacobian $\lambda = 4t_2^3 + t_1$ and $\eta = \pa/\pa t_2$. 
We follow the method of \cite{Ishikawa95}. 
Set 
$$
{\textstyle 
U =  t_2^4 + t_1t_2, \ 
U_1 = \frac{4}{5}t_2^5 + \frac{1}{2}t_1t_2^2, \ 
U_2 = \frac{2}{3}t_2^6 + \frac{1}{3}t_1t_2^3, \ 
U_3 = \frac{4}{7}t_2^7 + \frac{1}{4}t_1t_2^4. 
}
$$
The third component $f_3$ is written as 
$$
f_3 = A\circ g + (B\circ g)U_1 + (C\circ g)U_2 + (D\circ g)U_3. 
$$
Then the condition $\ord_a^\eta(f_3) = 5$ implies that $B(0, 0) \not= 0$. 
We may suppose $f = (g, f_3)$ with $f_3 = U_1 + \Phi, 
\Phi =  (B\circ g)U_1 + (C\circ g)U_2 + (D\circ g)U_3, B(0, 0) = 0$. Then 
the family $F_s = (g, \ U_1 + s\Phi)$ is trivialized by ${\mathcal A}$-equivalence. 
Thus $f = F_1$ is ${\mathcal A}$-equivalent to $F_0$ which is the normal form of (1'). 
Therefore (2) implies (1'). The converse is clear. Hence (1') and (2) are equivalent. 
\QED

\

As for openings of the lips $(t_1, t_2) \to (t_1, t_2^3 + t_1^2t_2)$ (see \cite{Saji}), we have 

\bet
{\rm (Recognition of cuspidal lips)} \ 
Let $(\R^2, a) \to (\R^3, b)$ be a frontal of corank $1$. 
Then the following conditions are equivalent to each other: 

{\rm (1)} $f$ is ${\mathcal A}$-equivalent to cuspidal lips i.e. 
$
{\textstyle 
(t_1, t_2) \to (t_1, t_2^3 + t_1^2t_2, \frac{3}{4}t_2^4 + \frac{1}{2}t_1^2t_2^2).
}
$ 

{\rm (2)} $f$ is a front and $\lambda$ is ${\mathcal K}$-equivalent to the germ 
$(t_1, t_2) \mapsto t_1^2 + t_2^2$ at the origin. 

{\rm (3)} $\lambda$ is ${\mathcal K}$-equivalent to the germ 
$(t_1, t_2) \mapsto t_1^2 + t_2^2$ at the origin, and $\ord_a^\eta(f_3) = 4$. 
\ent

\Proof
The equivalence of (1) and (2) is proved in \cite{IST}. 
Under the condition that 
$\lambda$ is ${\mathcal K}$-equivalent to the germ 
$(t_1, t_2) \mapsto t_1^2 + t_2^2$ at the origin, the condition $\ord_a^\eta(f_3) = 4$ is equivalent
to that the Legendre lift $\widetilde{f}$ is an immersion. Thus we have the equivalence of (2) and (3). 
\QED

\ber
{\rm 
Cuspidal lips never appear as singularities of tangent surfaces. 
}
\enr

\

We conclude the paper by imposing open questions: 

\

\noindent
{\bf Question 1.} When does ${\mathcal J}$-equivalence imply ${\mathcal A}$-equivalence ? 

\ber
{\rm 
For immersions, folds, cusps, lips, beaks, swallowtails : $(\R^2, 0) \to (\R^2, 0)$, ${\mathcal J}$-equivalence implies ${\mathcal A}$-equivalence. 
}
\enr

\bee
{\rm (\cite{Rieger}, \cite{Kabata})
Let $f, f' : (\R^2, 0) \to (\R^2, 0)$ be defined by 
$f(t_1, t_2) = (t_1,   t_1t_2 + t_2^5 + t_2^7)$ (butterfly) and 
$f'(t_1, t_2) = (t_1, t_1t_2 + t_2^5)$ (elder butterfly). 
Then $f$ is not ${\mathcal A}$-equivalent to $f'$ and their recognition by Taylor coefficients is obtained 
by Kabata\cite{Kabata}. 
On the other hand we observe, by using the theory of implicit OED of first order, 
that $f$ is ${\mathcal J}$-equivalent to $f'$ in fact. 
Therefore we see that it is absolutely impossible to recognize them just in terms of 
kernel field $\eta$ and Jacobian $\lambda$. 
}
\ene

\noindent
{\bf Question 2.} When does ${\mathcal J}$-equivalence imply ${\mathcal K}$-equivalence ?

\

It can be shown, for frontals of corank $1$, that 
${\mathcal J}$-equivalence implies ${\mathcal K}$-equivalence under a mild condition: 

\bel
Let $f : (\R^n, a) \to (\R^m, b)$ and $f' : (\R^n, a') \to (\R^{m'}, b')$ be map-germs of corank $1$. 
If $f$ and $f'$ are ${\mathcal J}$-equivalent and $f$ is ${\mathcal K}$-finite, 
then $f$ and $f'$ are ${\mathcal K}$-equivalent, i.e. 
$(f^*{\mathfrak m}_b){\mathcal E}_a$ is transformed to $({f'}^*{\mathfrak m}_{b'}){\mathcal E}_{a'}$ 
by a diffeomorphism-germ $\sigma : (\R^n, a) \to (\R^n, a')$. Here ${\mathfrak m}_b \subset 
{\mathcal E}_b$ is the maximal ideal. 
The condition that $f$ is ${\mathcal K}$-finite means that 
$\dim_{\R}({\mathcal E_a}/(f^*{\mathfrak m}_b){\mathcal E}_a) < \infty$. 
\enl

\Proof
By the assumption, $f$ is ${\mathcal A}$-equivalent to $g : (\R^n, 0) \to (\R^m, 0)$ 
of form 
$(t_1, \dots, t_{n-1}, \varphi_n(t), \dots, \varphi_m(t))$ for some $\varphi_i \in {\mathcal E}_0, n \leq i \leq m$. 
Then $g^*({\mathfrak m}_0){\mathcal E}_0$ is generated by $t_1, \dots, t_{n-1}, t_n^\ell$ 
for some $\ell$ and $\ell$ is uniquely determined by the minimum of 
orders of $\varphi_n(0, t_n), \dots, \varphi_m(0, t_n)$ for $t_n$ at $0$.  
On the other hand,  the Jacobi module ${\mathcal J}_g$ is generated by 
$dt_1, \dots, dt_{n-1}, (\pa \varphi_n/\pa t_n)d t_n, \dots, (\pa \varphi_n/\pa t_n)d t_n$, 
and the minimum of orders of $(\pa\varphi_n/\pa t_n)(0, t_n), \dots, (\pa\varphi_m/\pa t_n)(0, t_n)$ for $t_n$ at $0$ 
is invariant under ${\mathcal J}$-equivalence. Therefore ${\mathcal K}$-equivalence class is also invariant under 
${\mathcal J}$-equivalence. 
\QED


{\small 
\begin{thebibliography}{99}

\bibitem{Arnold}
V.I. Arnol'd, 
{\it Lagrangian manifold singularities, asymptotic rays and the open swallowtail, }
Funct. Anal. Appl., {\bf 15} (1981). 235--246. 

\bibitem{CI}
S. Chino, S. Izumiya, 
{\it 
Lightlike developables in Minkowski 3-space, 
}
Demonstratio Mathematica {\bf 43--2} (2010), 387--399. 

\bibitem{FSUY}
S. Fujimori, K. Saji, M. Umehara, K. Yamada, 
{\it Singularities of maximal surfaces, }
Math. Z. {\bf 259} (2008), 827--848. 

\bibitem{Givental}
A.B. Givental, 
{\it Whitney singularities of solutions of partial differential equations, } 
J. Geom. Phys. {\bf 15--4} 
(1995), 353--368. 

%\bibitem{GG}
%M. Golubitsky, V. Guillemin, 
%{\it Stable Mappings and Their Singularities, } 
%Graduate Texts in Mathematics 14, Springer-Verlag, 
%(1973). 

\bibitem{Ishikawa83}
G. Ishikawa, 
{\it 
Families of functions dominated by distributions of ${\cal C}$-classes of mappings, } 
Ann. Inst. Fourier, {\bf 33-2} (1983), 199--217. 

\bibitem{Ishikawa95}
G. Ishikawa, 
{\it Developable of a curve and its determinacy relatively to 
the osculation-type}, 
Quarterly J. Math., {\bf 46} (1995), 437--451. 

\bibitem{Ishikawa96} 
G. Ishikawa, 
{\it Symplectic and Lagrange stabilities of open Whitney umbrellas}, 
Invent. math., {\bf 126--2} (1996), 215--234. 

\bibitem{Ishikawa12}
G. Ishikawa, 
{\it Singularities of tangent varieties to curves and surfaces, }
Journal of Singularities, {\bf 6} (2012), 54--83. 

\bibitem{Ishikawa13}
G. Ishikawa, 
{\it Tangent varieties and openings of map-germs, } 
RIMS K\={o}ky\={u}roku Bessatsu, {\bf B38} (2013), 119--137. 

\bibitem{Ishikawa14}
G. Ishikawa, 
{\it Openings of differentiable map-germs and unfoldings, } 
Topics on Real and Complex Singularities, 
Proc. of the 4th Japanese-Australian Workshop (JARCS4), Kobe 2011, 
World Scientific (2014), pp. 87--113. 

\bibitem{Ishikawa16}
G. Ishikawa, 
{\it Singularities of frontals,} 
in \lq\lq Singularities in Generic Geometry'', 
Advanced Studies in Pure Mathematics {\bf 78}, Math. Soc. Japan., (2018), pp. 55--106. 
\ 
\verb+http://arxiv.org/abs/1609.00488/+

\bibitem{IJ}
G. Ishikawa, S. Janeczko, 
{\it 
Symplectic bifurcations of plane curves and isotropic liftings, } 
Quarterly J. of Math. Oxford, {\bf 54} (2003), 73--102. 


\bibitem{IM}
G. Ishikawa, Y. Machida, 
{\it Singularities of improper affine spheres and surfaces of constant Gaussian curvature, }
International Journal of Mathematics, {\bf 17--3} (2006), 269--293. 

\bibitem{IMT11}
G. Ishikawa, Y. Machida, M. Takahashi, 
{\it 
Asymmetry in singularities of tangent surfaces in contact-cone Legendre-null duality, 
}
Journal of Singularities, {\bf 3} (2011), 126--143. 

\bibitem{IMT16}
G. Ishikawa, Y. Machida, M. Takahashi, 
{\it Singularities of tangent surfaces in Cartan's split $G_2$-geometry, } 
Asian Journal of Mathematics, {\bf 20--2}, (2016), 353--382. 
%\\
%\verb+http://eprints3.math.sci.hokudai.ac.jp/2236/+

\bibitem{IMT18}
G. Ishikawa, Y. Machida, M. Takahashi, 
{\it 
Null frontal singular surfaces in Lorentzian 3-spaces, }
in preparation. 


\bibitem{IY17}
G. Ishikawa, T. Yamashita, 
{\it 
Singularities of tangent surfaces to generic space curves, } 
Journal of Geometry, {\bf 108} (2017), 301--318. 


\bibitem{IY18}
G. Ishikawa, T. Yamashita, 
{\it Singularities of tangent surfaces to directed curves, }
Topology and its Applications, {\bf 234} (2018), 198--208. 

%\bibitem{INS}
%S. Izumiya, T. Nagai, K. Saji, 
%{\it Great circular surfaces in the three-sphere, }
%Diff. Geom. and its Appl., {\bf 29} (2011) 409--425. 

\bibitem{IST}
S. Izumiya, K. Saji, M. Takahashi, 
{\it Horospherical flat surfaces in Hyperbolic 3-space, }
J. Math. Soc. Japan
{\bf 62-3} (2010), 789--849. 

\bibitem{Kabata}
Y. Kabata, 
{\it Recognition of plane-to-plane map-germs,}
Topology and its Applications, {\bf 202-1} (2016), 
216--238. 

\bibitem{KRSUY}
M. Kokubu, W. Rossman, K. Saji, M. Umehara, K. Yamada, 
{\it Singularities of flat fronts in hyperbolic space, }
Pacific J. of Math. {\bf 221-2} (2005), 303--351. 

\bibitem{Mather}
J.N. Mather, 
%{\it 
%Stability of $C^\infty$ mappings II: Infinitesimal stability implies stability, }
%Ann. of Math. {\bf 89} (1969), 254--291. 
{\it 
Stability of $C^\infty$ mappings III: 
Finitely determined map-germs,}
Publ. Math. I.H.E.S., {\bf 35} (1968), 279--308. 
%{\it 
%Stability of $C^\infty$ mappings IV: 
%Classification of stable germs by $\R$ algebras, 
%}
%Publ. Math. I.H.E.S., {\bf 37} (1969), 223--248. 
%{\it 
%Stability of $C^\infty$ mappings V: 
%Transversality, 
%}
%Advances in Math., {\bf 4} (1970), 301--336. 
%{\it 
%Stability of $C^\infty$ mappings VI: 
%The nice dimensions, 
%}
%Lecture Notes in Math., {\bf 192} (1971), Springer, Berlin, 
%pp. 207--253. 

\bibitem{Rieger}
J.H. Rieger, {\it 
Families of maps from the plane to the plane, }
 J. London Math. Soc., {\bf 36} (1987), 351--369. 


\bibitem{Saji}
K. Saji, 
{\it Criteria for singularities of smooth maps from the plane into the plane and their applications, }
Hiroshima Math. J. {\bf 40} (2010), 229--239. 


\bibitem{SUY}
K. Saji, M. Umehara, K. Yamada, 
{\it $A_k$ singularities of wave fronts, }
Math. Proc. Camb. Philos. Soc. {\bf 146-3} (2009), 731--746.


\bibitem{Sharp}
R.W. Sharp, 
{\it 
Differential Geometry: Cartan's Generalization of Klein's Erlangen Program}, 
 Graduate Texts in Mathematics {\bf 166}, Springer, (2000). 

\bibitem{Whitney}
H. Whitney, {\it On singularities of mappings of Euclidean spaces I,  
Mappings of the plane into the plane, } 
Ann. of Math. {\bf 62} (1955), 374--410. 

\end{thebibliography}

}
\end{document}